\documentclass[11pt,a4paper]{article}

\usepackage{epic,eepic,epsfig,amsmath,amssymb,amsthm}
\usepackage{t1enc}

\usepackage[francais,english]{babel}

\usepackage{caption}

\usepackage{color}

\usepackage{bbm}

\def\virgp{\raise 2pt\hbox{,}}

\renewcommand{\geq}{\geqslant}
\renewcommand{\leq}{\leqslant}
\def\N{{\mathbb N}}
\def\C{\mathbb{C}}
\def\R{{\mathbb R}}

\def\virgp{\raise 2pt\hbox{,}}
\def\cdotpv{\raise 2pt\hbox{;}}

\def\1{\mathbbm{1}}
\newcommand{\ds}{\displaystyle}

\newtheorem{theorem}{Theorem}[section]

\newtheorem{proposition}[theorem]{Proposition}

\newtheorem{pte}[theorem]{Property}

\theoremstyle{remark}
\newtheorem{remark}{Remark}[section]

\theoremstyle{definition}
\newtheorem{definition}{Definition}[section]

\newtheorem*{notation}{Notation}

\theoremstyle{definition}

\theoremstyle{definition}

\begin{document}
	
	\title{A spectral study of the Minkowski Curve}
	
	\author{Nizare Riane, Claire David}

	\maketitle
	\centerline{Sorbonne Universit\'e}
	
	\centerline{CNRS, Laboratoire Jacques-Louis Lions, 4, place Jussieu 75005, Paris, France}

	

	\maketitle
	\vskip 1cm

	\begin{abstract}
		In the following, we give an explicit construction of a Laplacian on the Minkowski curve, with energy forms that bear the geometric characteristic of the structure. The spectrum of the Laplacian is obtained by means of spectral decimation.
		
	\end{abstract}

	\vskip 1cm
	\noindent \textbf{Keywords}: Laplacian - Minkowski Curve - Einstein relation - Spectral decimation.

	\vskip 1cm
	
	\noindent \textbf{AMS Classification}:  37F20- 28A80-05C63.
	\vskip 1cm
	
	\vskip 1cm
	
	\section{Introduction}

	\hskip 0.5cm The so-called Minkowski curve, a fractal meandering one, whose origin seems to go back to Hermann Minkowski, but is nevertheless difficult to trace precisely, appears as an interesting fractal object, and a good candidate for whom aims at catching an overview of spectral properties of fractal curves.\\

	The curve is obtained through an iterative process, starting from a straight line which is replaced by eight segments, and, then, repeating this operation. One may note that the length of the curve grows faster than the one of the Koch one. What are the consequences in the case of a diffusion process on this curve ?\\
	
	The topic is of interest. One may also note that, in electromagnetism, fractal antenna, specifically, on the model of the aforementioned curve,  which miniaturized design turn out to perfectly fit wideband or broadband transmission, are increasingly used.\\

	It thus seemed interesting to us to build a specific Laplacian on those curves. To this purpose, the analytical approach initiated by~J.~Kigami~\cite{Kigami1989},~\cite{Kigami1993}, taken up, developed and popularized by R.~S.~Strichartz~\cite{Strichartz1999}, \cite{StrichartzLivre2006}, appeared as the best suited one. The Laplacian is obtained through a weak formulation, by means of Dirichlet forms, built by induction on a sequence of graphs that converges towards the considered domain. It is these Dirichlet forms that enable one to obtain energy forms on this domain. \\
	
	Laplacians on fractal curves are not that simple to implement. One must of course bear in mind that a fractal curve is topologically equivalent to a line segment. Thus, how can one make a distinction between the spectral properties of a curve, and a line segment ? Dirichlet forms solely depend on the topology of the domain, and not of its geometry. The solution is to consider energy forms more sophisticated than classical ones, by means of normalization constants that could, not only bear the topology, but, also, the geometric characteristics. One may refer to the works of U.~Mosco~\cite{UmbertoMosco2002} for the Sierpi\'{n}ski curve, and U.~Freiberg~\cite{UtaFreiberg2004}, where the authors build an energy form on non-self similar closed fractal curves, by integrating the Lagrangian on this curve. It is not the case of all existing works: in~\cite{UthayakumarDevi2014}, the authors just use topological normalization constants for the energy forms in stake.\\
	
	In the sequel, we give an explicit construction of a Laplacian on the Minkowski curve, with energy forms that bear the geometric characteristic of the structure. The spectrum of the Laplacian is obtained by means of spectral decimation, in the spirit of the works of M.~Fukushima and T.~Shima~ \cite{Fukushima1992}. On doing so, we choose three different methods. This enable us to initiate a detailed study of the spectrum.
	
	\begin{center}
		\includegraphics[scale=1]{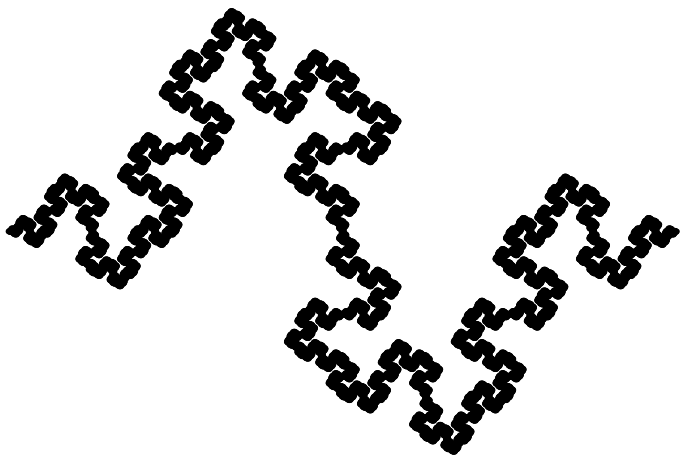}
		\captionof{figure}{Minkowski sausage.}
		\label{fig1}
	\end{center}

	\begin{figure}[h!]
		\center{\psfig{height=7cm,width=7cm,angle=0,file=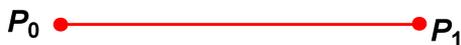}}\\
		\caption{The initial line segment.}
	\end{figure}

	\begin{figure}[h!]
		\center{\psfig{height=2cm,width=5cm,angle=0,file=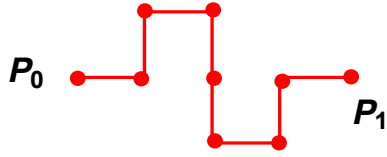}}\\
		\caption{The graph obtained at the first iteration.}
	\end{figure}

	\newpage
	\section{Framework of the study}

	In the sequel, we place ourselves in the Euclidean plane of dimension~2, referred to a direct orthonormal frame. The usual Cartesian coordinates are~$(x,y)$.\\

	\begin{notation}
		
		\noindent  Let us denote by~$P_0$ and ~$P_1$ the points:
		
		$$P_0= (0,0) \quad , \quad P_1= (1,0)$$
		
	\end{notation}
	
	\vskip 1cm

	\begin{notation}
		
		Let us denote by~$\theta \,\in\,]0,2\,\pi[$, $k>0$,~$T_1$, and~$T_2$ real numbers. Let us consider the rotation matrix:
		$${\cal R}_{O,\theta}=\left(
		\begin{matrix}
		\cos  \theta  & -   \sin  \theta \\
		\sin  \theta  &   \cos  \theta \\
		\end{matrix}
		\right)$$
		
		\noindent 	We introduce the iterated function system of the family of~maps from~$\R^2$ to~$\R^2 $:
		$$\left \lbrace f_{1},...,f_{8} \right \rbrace$$
		\noindent where, for any integer~$i$ belonging to~\mbox{$\left \lbrace 1,...,8 \right \rbrace$}, and any~$X \,\in\R^2$:
		$$f_i(X)=k\, {\cal R}_{O,\theta }
		X +
		\left(
		\begin{matrix}
		T_1\\
		T_2\\
		\end{matrix}
		\right)
		$$
		
	\end{notation}
	
	\vskip 1cm
	
	\begin{remark}
		\noindent If $0<k<1$, the family~$\left \lbrace f_{1},...,f_{8} \right \rbrace$ is a family of contractions from~$\R^2$ to~$\R^2$, the ratio of which is~$k$.
		
	\end{remark}
	
	\vskip 1cm

	\begin{pte}
		
		According to~\cite{Hutchinson1981}, there exists a unique subset $\mathfrak{MC} \subset \R^2$ such that:
		\[\mathfrak{MC} = \underset{  i=1}{\overset{8}{\bigcup}}\, f_i(\mathfrak{MC})\]
		\noindent which will be called the Minkowski Curve.\\
		
		\noindent One may note that this curve is not strictly self-similar in the sense defined by B.~Mandelbrot~\cite{Mandelbrot1977},~\cite{Hutchinson1981}:
		
		$$ \exists\, (i,j)\, \in~\left \lbrace f_{1},...,f_{8} \right \rbrace^2, \, i \neq j\, : \quad  f_i(\mathfrak{MC}) \cap f_j(\mathfrak{MC}) \neq \emptyset$$


	\end{pte}

	\vskip 1cm

	\begin{definition}
		
		\noindent We will denote by~$V_0$ the ordered set, of the points:
		
		$$\left \lbrace P_{0},P_{1}\right \rbrace$$

		\noindent The set of points~$V_0$, where, for any~$i$ of~\mbox{$\left \lbrace  0,1  \right \rbrace$}, the point~$P_0$ is linked to the point~$P_{1}$, constitutes an oriented graph, that we will denote by~$ {\mathfrak {MC}}_0$.~$V_0$ is called the set of vertices of the graph~$ {\mathfrak {MC}}_0$.\\
		
		\noindent For any strictly positive integer~$m$, we set:
		$$V_m =\underset{  i=1}{\overset{8}{\bigcup}}\, f_i  \left (V_{m-1}\right )$$

		\noindent The set of points~$V_m$, where the points of an $m$-cell are linked in the same way as ${\mathfrak {MC}}_0$, is an oriented graph, which we will denote by~$ {\mathfrak {MC}}_m$.~$V_m$ is called the set of vertices of the graph~$ {\mathfrak {MC}}_m$. We will denote, in the following, by~${\cal N}_m$ the number of vertices of the graph~$ {\mathfrak {MC}}_m$.

	\end{definition}
	
	\vskip 1cm
	
	\begin{notation}
		\noindent We introduce the following similarities, from~$\R^2$ to~$\R^2$, such that, for any~$X\,\in\,\R^2$:
		
		$$
		f_1(X) =\displaystyle \frac{1}{4} \,\left( {\cal R}_{O,0 }\,X +
		\left(
		\begin{matrix}
		0\\
		0\\
		\end{matrix}
		\right)
		\right) \quad , \quad
		f_2(X) =\displaystyle\frac{1}{4}\, \left( {\cal R}_{O,\frac{\pi}{2}}\,X +
		\left(
		\begin{matrix}
		1\\
		0\\
		\end{matrix}
		\right)
		\right)  \quad , \quad
		f_3(X) =\displaystyle\frac{1}{4} \,\left({\cal R}_{O,0}\,X +
		\left(
		\begin{matrix}
		1\\
		1\\
		\end{matrix}
		\right)
		\right) $$
		
		$$
		f_4(X) =\displaystyle\frac{1}{4} \,\left({\cal R}_{O,\frac{3\,\pi}{2}}\,X +
		\left(
		\begin{matrix}
		2\\
		1\\
		\end{matrix}
		\right)
		\right)  \quad , \quad
		f_5(X) =\displaystyle\frac{1}{4} \left({\cal R}_{O,\frac{3\,\pi}{2}}\,X +
		\left(
		\begin{matrix}
		2\\
		0\\
		\end{matrix}
		\right)
		\right)  \quad , \quad
		f_6(X) =\displaystyle\frac{1}{4} \,\left({\cal R}_{O,0}\,X +
		\left(
		\begin{matrix}
		2\\
		-1\\
		\end{matrix}
		\right)
		\right) $$
		
		$$
		f_7(X) =\displaystyle\frac{1}{4}\, \left({\cal R}_{O,\frac{\pi}{2}}\,X +
		\left(
		\begin{matrix}
		3\\
		-1\\
		\end{matrix}
		\right)
		\right)  \quad , \quad
		f_8(X) =\displaystyle\frac{1}{4}\, \left({\cal R}_{O,0}\,X +
		\left(
		\begin{matrix}
		3\\
		0\\
		\end{matrix}
		\right)
		\right) \\
		$$
	\end{notation}
	
	\vskip 1cm
	\begin{pte}
		
		The set of vertices~$\left (V_m \right)_{m \in\N}$ is dense in~$ {\mathfrak {MC}}$.
		
	\end{pte}

	\vskip 1cm
	\begin{proposition}
		Given a natural integer~$m$, we will denote by~$\mathcal{N}_m$ the number of vertices of the graph ${\mathfrak {MC}}_m$. One has:
		$$\mathcal{N}_0  =2$$
		
		\noindent and, for any strictly positive integer~$m$:
		$$\mathcal{N}_m   =8\times \mathcal{N}_{m-1}- 7$$
	\end{proposition}
	
	\vskip 1cm
	\begin{proof}
		The proposition results from the fact that, for any strictly positive integer~$m$, each graph~${\mathfrak {MC}}_m$ is the union of eight copies of ${\mathfrak {MC}}_{m-1}$: every copy shares one vertex  with its predecessor.\\
	\end{proof}
	
	\vskip 1cm

	\begin{definition}\textbf{Consecutive vertices of~${\mathfrak {MC}}$ }\\
		
		\noindent Let us set:~$P_2=P_0$.
		\noindent Two points~$X$ and~$Y$ of~${\mathfrak {MC}}$ will be called \textbf{\emph{consecutive vertices}} of~${\mathfrak {MC}} $ if there exists a natural integer~$m$, and an integer~$j $ of~\mbox{$\left \lbrace  0,1 \right \rbrace$}, such that:
		
		$$X = \left (f_{i_1}\circ \hdots \circ f_{i_m}\right)(P_j) \quad \text{and} \quad Y = \left (f_{i_1}\circ \hdots \circ f_{i_m}\right)(P_{j+1})
		\qquad \left ( i_1,\hdots, i_m \right ) \,\in\,\left \lbrace  1,\dots,8 \right \rbrace^m $$

	\end{definition}
	
	\vskip 1cm

	\begin{definition}\textbf{Edge relation, on the graph~$\mathfrak {MC}$}\\

		\noindent Given a natural integer~$m$, two points~$X$ and~$Y$ of~${\mathfrak {MC}}_m $ will be called \emph{\textbf{adjacent}} if and only if~$X$ and~$Y$ are two consecutive vertices of~${\mathfrak {MC}}_m $. We will write:
		
		$$X \underset{m }{\sim}  Y$$

	\end{definition}

\vskip 1cm

	\begin{definition}
		\noindent For any positive integer~$m$, the~$ {\cal N}_m$ consecutive vertices of the graph~$ {\mathfrak {MC}}_m  $ are, also, the vertices of~two triangles~${\cal T}_{m,left}$,~${\cal T}_{m,right}$, and~$2 \times 8^{m-1}$ squares~${\cal Q}_{m,j}$,~\mbox{$1 \leq j \leq 2 \times 8^{m-1}$}. For any integer~$j$ such that~\mbox{$1 \leq j \leq 2 \times 8^{m-1} $}, one obtains each square by linking the point number~$j$ to the point number~$j+1$ if~\mbox{$j = i \, \text{mod } 9$},~\mbox{$1 \leq i \leq   3$}, and the point number~$j$ to the point number~$j-3$ if~\mbox{$j =4 \, \text{mod } 9$}, and the triangles by linking the point number~$j$ to the point number $j+1$ and $j+2$ if $j=0 \,\text{mod } 9$. One has to consider those polygons as semi-closed ones, since, for any of those~$4-$gons, the starting vertex, i.e. the point number~$j$, is not connected, on the graph~$ {\mathfrak {MC}}_m  $, to the extreme one, i.e. the point number~$j-3$. In the same way, the extreme vertices of the~two triangles~${\cal T}_{m,left}$,~${\cal T}_{m,right}$, are not connected, on the graph~$ {\mathfrak {MC}}_m  $. \\
		\noindent The reunion of those triangles and squares generate a Borel set of~$\R^2$.\\

		\begin{figure}[h!]
			\center{\psfig{height=5cm,width=6.5cm,angle=0,file=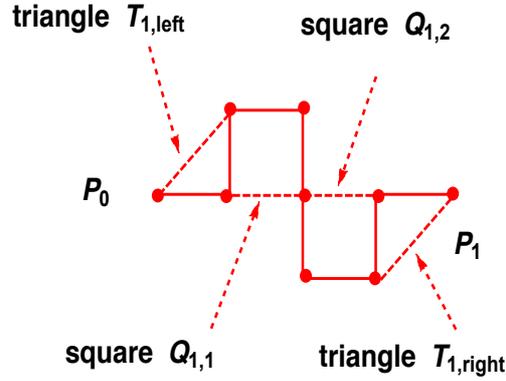}}\\
			\caption{The triangles~${\cal T}_{1,left}$,~${\cal T}_{1,right}$, and the squares~${\cal Q}_{1,1}$,~${\cal Q}_{1,2}$.}
		\end{figure}
		
	\end{definition}

\vskip 1cm

	\begin{definition}\textbf{Quasi-quadrilateral domain delimited by the graph~$ {\mathfrak {MC}}_m  $,~$m\,\in\,\N $}\\
		
		\noindent For any natural integer~$m$, we will call \textbf{quasi-quadrilateral domain delimited by~$ {\mathfrak {MC}}_m  $},
		and denote by~\mbox{$ {\cal D} \left ( {{\mathfrak {MC}}_m}\right) $}, the reunion of the triangles~${\cal T}_{m,j}$, and the squares~${\cal Q}_{m,j}$.
	\end{definition}
	\vskip 1cm

	\begin{definition}\textbf{Quasi-quadrilateral  domain delimited by the Minkowski Curve~$  {\mathfrak {MC}}$ }\\
		
		\noindent We will call \textbf{quasi-quadrilateral  domain delimited by~$  {\mathfrak {MC}}$}, and denote by~\mbox{$ {\cal D} \left ( {\mathfrak {MC}}\right) $}, the limit:
		$$ {\cal D} \left ( {\mathfrak {MC}}\right)  = \displaystyle \lim_{n \to + \infty} {\cal D} \left ({\mathfrak {MC}}_m \right) $$

		\noindent which has to be understood in the following way: given a continuous function~$u$ on the curve~$\mathfrak {MC}$, and a measure with full support~$\mu$ on~$\R^2$, then:
		\footnotesize
		$$\displaystyle \int_{ {\cal D} \left ( {\mathfrak {MC}}\right) } \! u\,d\mu  =
		\displaystyle \lim_{m \to + \infty}
		\left \lbrace \displaystyle \frac{1}{3}  \displaystyle \sum_{X \, \text{vertex of }{\cal T}_{m,left} } \! \!  \!\! \! \! \! \! u\left ( X \right) \,\mu \left (  {\cal T}_{m,left}  \right)
		+\displaystyle \frac{1}{4} 	\displaystyle \sum_{j=1}^{2 \times 8^{m-1}  } \! \!  \displaystyle \sum_{X \, \text{vertex of }{\cal Q}_{m,j} } \! \! \! \!  \! \! \!\! u\left ( X \right) \,\mu \left (  {\cal Q}_{m,j}  \right) 
		+\displaystyle \frac{1}{3} \displaystyle \sum_{X \, \text{vertex of }{\cal T}_{m,right} }\! \!  \!\! \! \! \! \! u\left ( X \right) \,\mu \left (  {\cal T}_{m,right}  \right)  \right \rbrace $$
		\normalsize
		
	\end{definition}
	
	\vskip 1cm

	\begin{remark}
		
		One may wonder why consider such a quasi-quadrilateral domain, whereas we presently deal with a curve ? If it seemed, to us, the most natural way to perform integration on our fractal curve, there are, also, additionnal reasons. Part of an explanation can be found in the paper of~U.~Mosco~\cite{UmbertoMosco2002}, where the author suggests, in the case of the Sierpi\'nski curve,  to "fill each small simplex (...) not only with its edges, but with the whole portion of
		the limit curve which it encompasses." It is also natural, in so far as it enables one to take into account the self-similar patterns that appear thanks to the curve, patterns, the measure of which plays the part of a pound. This joins the seminal work of J.~Harrison~et al.~\cite{JennyHarrison1991}, ~\cite{JennyHarrison1992}, ~\cite{JennyHarrison1993}.

	\end{remark}
	\vskip 1cm
	\begin{definition}\textbf{Word, on the Minkowski curve ~${\mathfrak {MC}}$}\\

		\noindent Let~$m  $ be a strictly positive integer. We will call \textbf{number-letter} any integer~${\cal W}_i$ of~\mbox{$\left \lbrace 1,\dots,8 \right \rbrace $}, and \textbf{word of length~$|{\cal W}|=m$}, on the graph~$\mathfrak {MC}$, any set of number-letters of the form:

		$${\cal W}=\left ( {\cal W}_1, \hdots, {\cal W}_m\right)$$
		
		\noindent We will write:
		
		$$f_{\cal W}= f_{{\cal W}_1} \circ \hdots \circ  f_{{\cal W}_m}  $$
		
	\end{definition}
	
	\vskip 1cm

	\begin{proposition}\textbf{Adresses, on the the Minkowski Curve}\\

		\begin{itemize}
			\item[i.] Every~$P_i$,~$i\,\in\, \{0,1\}$, has exactly one neighbor.
			\item[ii.] Given a strictly positive integer~$m$, every~$X \,\in \,V_m\setminus V_0$ has exactly two neighbors and two addresses 
			$$X=f_{\left ( {\cal W}_1, \hdots, {\cal W}_m\right)}(P_0)=f_{\left ( {\cal W'}_1, \hdots, {\cal W'}_m\right)}(P_1)$$
			
			where~${\cal W}=\left ( {\cal W}_1, \hdots, {\cal W}_m\right)$ and ${\cal W'}=\left ( {\cal W'}_1, \hdots, {\cal W'}_m\right)$ denote two appropriate words of length~$m $, on the graph~${\mathfrak {MC}}_m  $~(see~\cite{StrichartzLivre2006}).

		\end{itemize}
		
	\end{proposition}
	
	\vskip 1cm
	
	\begin{proposition}
		Let us set:
		
		$$\ds{V_\star =\underset{{m\in \N}}\bigcup \,V_m} $$
		
		\noindent The set $V_{\star}$ is dense in~$\mathfrak{MC}$.\\
	\end{proposition}
	
	\newpage
	
	\section{Energy forms, on the Minkowski Curve}
	
	\subsection{Dirichlet forms}

	\begin{definition}\textbf{Dirichlet form, on a finite set} (We refer to~\cite{Kigami2003})\\
		
		\noindent Let~$V$ denote a finite set~$V$, equipped with the usual inner product which, to any pair~$(u,v)$ of functions defined on~$V$, associates:
		
		$$(u,v)= \displaystyle \sum_{P\in  V} u(P)\,v(P)$$
		
		\noindent A \emph{\textbf{Dirichlet form}} on~$V$ is a symmetric bilinear form~${\cal E}$, such that:\\
		
		\begin{enumerate}
			
			\item For any real valued function~$u$ defined on~$V$:  ${\cal E}(u,u) \geq 0$.
			
			\item   $  {\cal {E}}(u,u)= 0$ if and only if~$u$ is constant on~$V$.

			\item For any real-valued function~$u$ defined on~$V$, if:
			$$ u_\star = \min\, (\max(u, 0) , 1)  $$
			
			\noindent i.e. :
			
			$$\forall \,p \,\in\,V \, : \quad u_\star(p)= \left \lbrace \begin{array}{ccc} 1 & \text{if}& u(p) \geq 1 \\u(p) & \text{si}& 0 <u(p) < 1 \\0  & \text{if}& u(p) \leq 0 \end{array} \right.$$
			
			\noindent then: ${ \cal{E}}(u_\star,u_\star)\leq { \cal{E}}(u,u)$ (Markov property).

		\end{enumerate}
		
	\end{definition}
	
	\vskip 1cm

	Let us now consider the problem of energy forms on our curve. Such a problem was studied by~U.~Mosco~\cite{UmbertoMosco2002}, who suggested to generalize Riemaniann models to fractals and relate the fractal analogous of gradient forms, i.e. the Dirichlet forms, to a metric that could reflect the fractal properties of the considered structure. The link is to be made by means of specific energy forms.\\
	
	There are two major features that enable one to characterize fractal structures:
	\begin{enumerate}
		\item[\emph{i}.] Their topology, i.e. their ramification.
		\item[\emph{ii}.] Their geometry.\\
		
	\end{enumerate}
	
	The topology can be taken into account by means of classical energy forms (we refer to~\cite{Kigami1989}, ~\cite{Kigami1993},~\cite{Strichartz1999}, \cite{StrichartzLivre2006}).\\
	As for the geometry, again, things are not that simple to handle.~U.~Mosco introduces a strictly positive parameter,~$\delta$, which is supposed to reflect the way ramification - or the iterative process that gives birth to the sequence of graphs that approximate the structure - affects the initial geometry of the structure. For instance, if $m$ is a natural integer,~$X$ and~$Y$ two points of the initial graph~$V_1$, and~$\cal M$ a word of length~$m$, the Euclidean distance~$d_{\R^2}(X,Y)$ between ~$X$ and~$Y$ is changed into the effective distance:
	
	$$\left ( d_{\R^2}(X,Y)\right )^\delta  $$
	
	This parameter~$\delta$ appears to be the one that can be obtained when building the effective resistance metric of a fractal structure~(see~\cite{StrichartzLivre2006}), which is obtained by means of energy forms. To avoid turning into circles, this means:
	\begin{enumerate}
		\item[\emph{i}.] either working, in a first time, with a value~$\delta_0$ equal to one, and, then, adjusting it when building the effective resistance metric ;
		\item[\emph{ii}.] using existing results, as done in~\cite{UtaFreiberg2004}.
		
	\end{enumerate}
	
	\noindent One may note that in a very interesting and useful remark, U.~Mosco puts the light on the relation that exists between the walk dimension~$D_W$ of a self-similar set, and~$\delta$:
	
	$$D_W=2 \,\delta$$
	
	\noindent which will enable us to obtain the required value of the constant~$\delta$.

	\vskip 1cm

	\begin{definition}\textbf{Walk dimension}\\
		
		\noindent Given a strictly positive integer~$N$, the dimension related to a random walk on a self-similar set with respect to~$N$ similarities, the ratio of which is equal to~$k\,\in\,]0,1[$, is given by:
		$$D_{W}=-\frac{\ln \mathbf{E}_A (\tau) }{\ln k} $$
		\noindent where $\mathbf{E}_A (\tau)$ is the mean crossing time of a random walk starting from a vertex~$A$ of the self-similar set.
	\end{definition}
	
	\vskip 1cm

	\begin{notation}
		
		In the sequel, we will denote by~$D_{W}\left(\mathfrak {MC}\right)$ the Walk dimension of the Minkowski Curve.

	\end{notation}

	\vskip 1cm
	
	\begin{remark}\textbf{Explicit computation of the Walk dimension of the Minkowski Curve}\\
		
		\noindent In the sequel, we follow the algorithm described in~\cite{UtaFreiberg2013}, which uses the theory of Markov chains (we refer to~\cite{KemenySnellMarkovChains}).\\
		
		\noindent To this purpose, we define $\mathbb{E}_X(\tau)$ to be the mean number of steps a simple random walk needs to reach a vertex~${B}\in V_0$ when starting at $X$. We consider the graph $\mathfrak{MI}_1$ and denote by $X_i$ for $i\in \{1,...,7\}$ the set of vertices of $V_1\setminus V_0$:
		
		\begin{center}
			\includegraphics[scale=1]{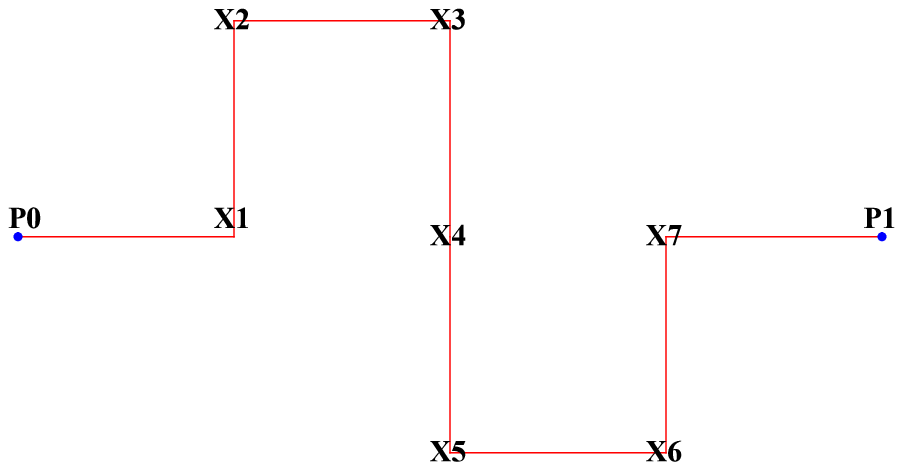}
			\captionof{figure}{The graph~$\mathfrak{MI}_1$}
			\label{fig5}
		\end{center}
		
		\noindent First, one has to introduce the adjacency matrix of the graph~$\mathfrak{MI}_1$:
		\footnotesize
		$$A_{\mathfrak{MI}_1}=
		\left ( \begin{array}{ccccccccc}
		0&1&0&0&0&0&0&0&0\\
		1&0&1&0&0&0&0&0&0\\
		0&1&0&1&0&0&0&0&0\\
		0&0&1&0&1&0&0&0&0\\
		0&0&0&1&0&1&0&0&0 \\
		0&0&0&0&1&0&1&0&0  \\
		0&0&0&0&0&1&0&1&0   \\
		0&0&0&0&0&0&1&0&1   \\
		0&0&0&0&0&0&0&1&0\\
		\end{array}
		\right)$$
		\normalsize
		\noindent One then build the related stochastic matrix:
		\footnotesize
		$$A_{\mathfrak{MI}_1}^{stoch}=
		\left ( \begin{array}{ccccccccc}
		0&1&0&0&0&0&0&0&0\\
		\frac{1}{2}&0&\frac{1}{2}&0&0&0&0&0&0\\
		0&\frac{1}{2}&0&\frac{1}{2}&0&0&0&0&0\\
		0&0 &\frac{1}{2}&0&\frac{1}{2}&0&0&0&0\\
		0&0&0 &\frac{1}{2}&0&\frac{1}{2}&0&0&0 \\
		0&0&0&0 &\frac{1}{2}&0&1&0&0  \\
		0&0&0&0&0 &\frac{1}{2}&0&\frac{1}{2}&0   \\
		0&0&0&0&0&0 &\frac{1}{2}&0&\frac{1}{2}   \\
		0&0&0&0&0&0&0&1&0 \\
		\end{array}
		\right)$$
		
		\normalsize
		
		\noindent By suppressing the lines and columns associated to one may call "cemetary states" (i.e. the fixed points of the contractions), one obtains the matrix:
		
		$$ M_{\mathfrak{MI}_1}=
		\left ( \begin{array}{ccccccccc}
		0&1&0&0&0&0&0&0\\
		\frac{1}{2}&0&\frac{1}{2}&0&0&0&0&0\\
		0&\frac{1}{2}&0&\frac{1}{2}&0&0&0&0\\
		0&0 &\frac{1}{2}&0&\frac{1}{2}&0&0&0\\
		0&0&0 &\frac{1}{2}&0&\frac{1}{2}&0&0 \\
		0&0&0&0 &\frac{1}{2}&0&1&0  \\
		0&0&0&0&0 &\frac{1}{2}&0&\frac{1}{2}  \\
		0&0&0&0&0&0 &\frac{1}{2}&0   \\
		\end{array}
		\right)$$
		
		\noindent Starting at $P_0$, one has:
		
		$$
		\mathbb{E}_{P_0}(\tau) =\mathbf{E}_{X_1}(\tau) +1 \quad , \quad
		\mathbb{E}_{P_1}(\tau) =0$$
		
		$$
		\mathbb{E}_{X_1}(\tau) =\displaystyle \frac{1}{2}\,\left (\mathbb{E}_{P_0}+\mathbf{E}_{X_2}(\tau)\right )+1\quad , \quad
		\mathbb{E}_{X_2}(\tau) =\displaystyle\frac{1}{2}\,\left (\mathbb{E}_{X_1}+\mathbb{E}_{X_3}(\tau)\right)+1\quad , \quad
		\mathbb{E}_{X_3}(\tau) =\displaystyle\frac{1}{2}\,\left (\mathbb{E}_{X_2}+\mathbb{E}_{X_4}(\tau)\right)+1$$
		
		$$\mathbb{E}_{X_4}(\tau) =\displaystyle\frac{1}{2}\,\left (\mathbb{E}_{X_3}+\mathbb{E}_{X_5}(\tau)\right)+1 \quad , \quad
		\mathbb{E}_{X_5}(\tau) =\displaystyle\frac{1}{2}\,\left (\mathbb{E}_{X_4}+\mathbb{E}_{X_6}(\tau)\right)+1\quad , \quad
		\mathbb{E}_{X_6}(\tau) =\displaystyle\frac{1}{2}\,\left (\mathbb{E}_{X_5}+\mathbb{E}_{X_7}(\tau)\right)+1$$
		$$\mathbb{E}_{X_7}(\tau) =\displaystyle\frac{1}{2}\,\left (\mathbb{E}_{P_1}+\mathbb{E}_{X_6}(\tau)\right)+1
		$$
		
		\noindent Let us introduce the vector of expected crossing-times:
		
		$$\mathbb{T}=
		\left(
		\begin{array}{c}
		\mathbb{E}_{P_0}(\tau)   \\
		\mathbb{E}_{X_1}(\tau) \\
		\mathbb{E}_{X_2}(\tau)  \\
		\mathbb{E}_{X_3}(\tau) \\
		\mathbb{E}_{X_4}(\tau) \\
		\mathbb{E}_{X_5}(\tau) \\
		\mathbb{E}_{X_6}(\tau) \\
		\mathbb{E}_{X_7}(\tau) \\
		\end{array}
		\right)$$
		
		From the theory of Markov chains, we have :
		
		$$\mathbb{T}=\left ( {\mathbb{E}}_{\mathfrak{MI}_1}- M_{\mathfrak{MI}_1}\right)^{-1} \mathbf{1}$$
		
		where $\mathbf{1}:=(1,...,1)^T$. By solving the above system, one gets:
		
		$$\mathbb{T} =
		\left(
		\begin{array}{c}
		64 \\
		63 \\
		60 \\
		55 \\
		48 \\
		39 \\
		28 \\
		15 \\
		\end{array}
		\right)
		$$
		
		\noindent which leads to:
		
		$$\mathbb{E}_{P_1}(\tau)=\displaystyle 64$$
		
		\noindent and:
		
		\vskip 1cm
		
		$$D_{W}\left(\mathfrak {MC}\right)=\displaystyle\frac{\ln 64 }{\ln 4 }=3$$

		\noindent one has thus:

		$$\delta= \displaystyle \frac{D_{W}\left(\mathfrak {MC}\right)}{2}= \displaystyle \frac{3}{2}$$

	\end{remark}
	
	\vskip 1cm

	\begin{definition}\textbf{Energy, on the graph~${\mathfrak {MC}}_m $,~$m \,\in\,\N$, of a pair of functions}\\

		\noindent Let~$m$ be a natural integer, and~$u$ and~$v$ two real valued functions, defined on the set
		
		$$V_m = \left \lbrace    X_1^m, \hdots,  X_{{\cal N}_m^{\cal S} }^m \right \rbrace $$
		
		\noindent of the~${\cal N}_m $ vertices of~${\mathfrak {MC}}_m   $.\\
		
		\noindent We introduce \textbf{the energy, on the graph~${\mathfrak {MC}}_m  $, of the pair of functions~$(u,v)$}, as:
		
		$$\begin{array}{ccc}
		{\cal{E}}_{{\mathfrak {MC}}_m  }(u,v)
		&= &  \displaystyle \sum_{i=1}^{{\cal N}_m -1}  \left (\displaystyle \frac{u \left (X_i^m \right)-u \left (X_{i+1}^m \right)}{d_{\R^2}^{\delta}(X,Y)}\right )\,
		\left (\displaystyle \frac{v \left (X_{i }^m \right)-v \left (X_{i+1}^m \right)}{ d_{\R^2}^{\delta}(X,Y)} \right )\\
		&= &  \displaystyle \sum_{i=1}^{{\cal N}_m -1} 4^{2\,m \,\delta}\,\left ( u \left (X_i^m \right)-u \left (X_{i+1}^m \right) \right )\,
		\left ( v \left (X_{i }^m \right)-v \left (X_{i+1}^m \right)  \right )\\
		\end{array}
		$$
		
		\noindent For the sake of simplicity, we will write it under the form:
		
		$$ {\cal{E}}_{{\mathfrak {MC}}_m  }(u,v)= \displaystyle \sum_{X  \underset{m }{\sim}  Y} 4^{2\, m \,\delta}\,\left (u(X)-u(Y)\right )\,\left(v(X)-v(Y)\right) $$

	\end{definition}
	
	\vskip 1cm

	\begin{proposition}\textbf{Harmonic extension of a function, on the Minkowski Curve - Ramification constant}\\
		
		\noindent For any integer~$m>1$, if~$u$ is a real-valued function defined on~$V_{m-1}$, its \textbf{harmonic extension}, denoted by~$ \tilde{u}$, is obtained as the extension of~$u$ to~$V_m$ which minimizes the energy:

		$$  {\cal{E}}_{{\mathfrak{MC}}_m }(\tilde{u},\tilde{u})= \displaystyle \sum_{X \underset{m }{\sim} Y} 4^{2\, m \,\delta}\,(\tilde{u}(X)-\tilde{u}(Y))^2 $$
		
		\noindent The link between~$   {\cal{E}}_{{\mathfrak {MC}}_{m } }$ and~$  {\cal{E}}_{{\mathfrak {MC}}_{m-1 } }$ is obtained through the introduction of two strictly positive constants~$r_m$ and~$r_{m-1}$ such that:

		$$   r_{m }\, \displaystyle \sum_{X \underset{m+1  }{\sim} Y} 4^{2\,( m+1) \,\delta}\, (\tilde{u}(X)-\tilde{u}(Y))^2 =  r_{m-1}\,4^{ 2\,(m -1)\,\delta}\, \sum_{X \underset{m-1 }{\sim} Y} (u(X)-u(Y))^2$$

		\noindent For the sake of simplicity, we will fix the value of the initial constant:~$r_0=1$.

		$$r = \displaystyle \frac{1}{r_{1 }} $$

		\noindent By induction, one gets:
		
		$$r_m=r_1^m =r^{-m} $$

		\noindent and:
		
		$$  {\cal{E}}_{m}(u)= r^{-m}\,4^{2\, m \,\delta}\, \sum_{X \underset{m }{\sim} Y}   (\tilde{u}(X)-\tilde{u}(Y))^2 $$

		\noindent Since the determination of the harmonic extension of a function appears to be a local problem, on the graph~${\mathfrak {MC}}_{ m-1} $, which is linked to the
		graph~$ {\mathfrak {MC}}_m  $ by a similar process as the one that links~${\mathfrak {MC}}_1$ to~$ {\mathfrak {MC}}_0$, one deduces, for any integer~$m>1$:

		$$ {\cal{E}}_{{\mathfrak {MC}}_m}(\tilde{u},\tilde{u})= r^{-m}\, 4^{2\, m \,\delta} \,  {\cal{E}}_{{\mathfrak {MC}}_{m-1}}(\tilde{u},\tilde{u})$$

		\noindent If~$v$ is a real-valued function, defined on~$V_{m-1}$, of harmonic extension~$ \tilde{v}$, we will write:
		
		$$  {\cal{E}}_{m}(u,v)= r^{-m}\,\displaystyle  \sum_{X \underset{m }{\sim} Y}  4^{2\, m \,\delta}\,(\tilde{u}(X)-\tilde{u}(Y)) \, (\tilde{v}(X)-\tilde{v}(Y))
		= r^{-m}\,\displaystyle  \sum_{X \underset{m }{\sim} Y}  4^{3\, m  }\,(\tilde{u}(X)-\tilde{u}(Y)) \, (\tilde{v}(X)-\tilde{v}(Y))$$

		\noindent The constant~$r^{-1}$, which can be interpreted as a topological one, will be called \textbf{ramification constant}.\\
		\noindent For further precision on the construction and existence of harmonic extensions, we refer to~\emph{\cite{Sabot1987}}.
	\end{proposition}

	\vskip 1cm
	
	\begin{definition}\textbf{Energy scaling factor}\\
		
		\noindent By definition, the \textbf{energy scaling factor} is the strictly positive constant~$\rho$ such that, for any integer~$m>1$, and any real-valued function~$u$ defined on~$V_{m }$:
		
		$$  {\cal{E}}_{{\mathfrak{MC}}_m }(u,u)= \rho\, {\cal{E}}_{{\mathfrak{MC}}_m }\left (u_{\mid V_{m-1}},u_{\mid V_{m-1}} \right)$$

	\end{definition}
	\vskip 1cm

	\begin{proposition}
		The energy scaling factor~$\rho$ is linked to the topology and the geometry of the fractal curve by means of the relation:
		
		$$\rho =\displaystyle \frac{4^{2\,\delta}}{8}=8 $$

		\noindent  One may note that a more general calculation of the energy scaling factor can be found in the work by~R.~Capitanelli~\cite{Capitanelli2002}.
		
	\end{proposition}
	
	\vskip 1cm

	\subsection{Explicit computation of the ramification constant}

	\subsubsection{Direct method}
	Let us denote by~$u$ a real-valued, continuous function defined on~$V_0=\{P_0,P_1\}$, and by~$\tilde{u}$ its harmonic extension to~$V_1$. Let us set: $u(P_0)=A$, $u(P_1)=B$. We recall that the energy on $V_0$ is given by :
	$$
	E_0(u)=(A-B)^2 $$
	\noindent  We will denote by $U_{1}$, $U_{2}$, $U_{3}$, $U_{4}$, $U_{5}$, $U_{6}$, $U_{7}$, $U_{16}$, the respective images, by~$\tilde{u}$, of the vertices of $V_1\setminus V_0$.\\
	
	\noindent One has then:
	\footnotesize
	$$
	E_1(\tilde{u}) =(A-U_1)^2+(U_1-U_2)^2+(U_2-U_3)^2+(U_3-U_4)^2+(U_4-U_5)^2+(U_5-U_6)^2 +(U_6-U_7)^2+(U_7-B)^2$$
	
	\normalsize

	\noindent The minimum of this quantity is such that:
	
	\[ \mathbf{U}=\mathbf{A}^{-1}\,\mathbf{b}\]

	\noindent where the matrix~$A$ is given by:
	
	\[
	\mathbf{A} =
	\left(
	\begin{array}{ccccccc}
	2 & -1 & 0 & 0 & 0 & 0 & 0 \\
	-1 & 2 & -1 & 0 & 0 & 0 & 0 \\
	0 & -1 & 2 & -1 & 0 & 0 & 0 \\
	0 & 0 & -1 & 2 & -1 & 0 & 0 \\
	0 & 0 & 0 & -1 & 2 & -1 & 0 \\
	0 & 0 & 0 & 0 & -1 & 2 & -1 \\
	0 & 0 & 0 & 0 & 0 & -1 & 2 \\
	\end{array}
	\right)
	\]
	
	\noindent the vectors~$\mathbf{U}$ and~$\mathbf{b}$ by:
	$$
	\mathbf{U} =
	\begin{pmatrix}
	U_{1} \\
	U_{2} \\
	U_{3} \\
	U_{4} \\
	U_{5} \\
	U_{6} \\
	U_{7} \\
	\end{pmatrix}
	\quad , \quad
	\mathbf{b} =
	\left(
	\begin{array}{c}
	A \\
	0 \\
	0 \\
	0 \\
	0 \\
	0 \\
	B \\
	\end{array}
	\right)$$
	
	\noindent One obtains:
	
	\[
	\mathbf{U}=
	\left(
	\begin{array}{c}
	\frac{7 A}{8}+\frac{B}{8} \\
	\frac{3 A}{4}+\frac{B}{4} \\
	\frac{5 A}{8}+\frac{3 B}{8} \\
	\frac{A}{2}+\frac{B}{2} \\
	\frac{3 A}{8}+\frac{5 B}{8} \\
	\frac{A}{4}+\frac{3 B}{4} \\
	\frac{A}{8}+\frac{7 B}{8} \\
	\end{array}
	\right)
	\]
	
	\noindent By substituting these values in the expression of the energy expression, one obtains:
	$$
	E_1(\tilde{u})= \displaystyle \frac{(A-B)^2}{8}=\displaystyle \frac{1}{8}\, E_0(u)
	$$
	
	\noindent Thus:
	
	$$r=\displaystyle \frac{1}{8}$$

	\vskip 1cm

	\subsubsection{A second method, using Einstein's relation}

	\begin{definition}\textbf{Hausdorff dimension of a self-similar set with respect to~$N\,\in\,\N^\star$ similarities}\\
		
		\noindent Given a strictly positive integer~$N$, the Hausdorff dimension of a self-similar set with respect to~$N$ similarities, the ratio of which is equal to~$k\,\in\,]0,1[$, is given by:
		$$D_{H}=-\frac{\ln N}{\ln k} $$
		
	\end{definition}

	\vskip 1cm

	\begin{notation}
		
		In the sequel, we will denote by
		$$D_{H}\left(\mathfrak {MC}\right)=-\frac{\ln 8}{\ln k}  $$
		
		\noindent  the Hausdorff dimension of the Minkowski Curve.

	\end{notation}

	\vskip 1cm

	\begin{definition}\textbf{Spectral dimension of a self-similar set with respect to~$N\,\in\,\N^\star$ similarities}\\
		
		\noindent Given a strictly positive integer~$N$, the spectral dimension of a self-similar set with respect to~$N$ similarities, the energy scaling factor of which is equal to $\rho$, is given by:
		$$D_{S}=\frac{2\, \ln N}{\ln(N\times \rho)}  $$

	\end{definition}

	\vskip 1cm

	\begin{notation}
		
		In the sequel, we will denote by
		$$D_{S}\left(\mathfrak {MC}\right) = \frac{6\, \ln 2}{\ln(8\times \rho)}$$
		
		\noindent  the spectral dimension of the Minkowski Curve.

	\end{notation}

	\vskip 1cm
	
	\begin{theorem}\textbf{Einstein relation}\\
		
		\noindent Given a strictly positive integer~$N$, and a self-similar set with respect to~$N$ similarities, the ratio of which is equal to~$k\,\in\,]0,1[$, one has the so-called Einstein relation  between the walk dimension, the Hausdorff dimension, and the spectral dimension of the set:
		$$  D_{H}=\frac{D_{S}\, D_{W}}{2} $$
	\end{theorem}
	
	\vskip 1cm

	\begin{pte}
		
		\noindent Given a strictly positive integer~$N$, and a self-similar set with respect to~$N$ similarities, the ratio of which is equal to~$k\,\in\,]0,1[$, the energy scaling factor~$\rho$, which solely depends on the topology, and, therefore, not of the value of the contraction ratio, is given by:
		$$ \rho=\frac{\mathbb{E}_A(T)}{N} $$
	\end{pte}
	
	\begin{center}
		\includegraphics[scale=1]{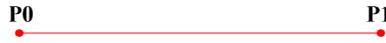}
		\captionof{figure}{The graph~$\mathfrak{MC}_0$}
		\label{fig4}
	\end{center}

	\vskip 1cm

	\subsection{Specific Dirichlet and energy forms, for the Minkowski Curve}
	
	\begin{definition}\textbf{Dirichlet form, for a pair of continuous functions defined on the Minkowski Curve~$\mathfrak {MC}$}\\
		
		\noindent We define the Dirichlet form~$\cal{E}$ which, to any pair of real-valued, continuous functions~$(u,v)$ defined on~$\mathfrak {MC}$, associates, subject to its existence:
		
		$$
		{\cal{E}} (u,v)= \displaystyle \lim_{m \to + \infty} {\cal{E}}_{m }\left (u_{\mid V_m},v_{\mid V_m}\right)=
		\displaystyle \lim_{m \to + \infty}\displaystyle \sum_{X  \underset{m }{\sim}  Y} r^{-m}\,4^{ 2\,m \,\delta}\, \left (u_{\mid V_m}(X)-u_{\mid V_m}(Y)\right )\,\left(v_{\mid V_m}(X)-v_{\mid V_m}(Y)\right) $$

	\end{definition}

	\vskip 1cm
	\begin{definition}\textbf{Normalized energy, for a continuous function~$u$, defined on~$\mathfrak {MC} $}\\
		\noindent Taking into account that the sequence~$\left (\mathcal{E}_m\left ( u_{\mid V_m} \right)\right)_{m\in\N}$ is defined on
		$$\ds{V_\star =\underset{{i\in \N}}\bigcup \,V_i}$$
		
		\noindent one defines the normalized energy, for a continuous function~$u$, defined on~$\mathfrak {MC}$, by:

		$$\mathcal{E}(u)=\underset{m\rightarrow +\infty}\lim \mathcal{E}_m \left ( u_{\mid V_m} \right)$$

	\end{definition}
	
	\vskip 1cm

	\begin{pte}
		\noindent The Dirichlet form~$\cal{E}$ which, to any pair of real-valued, continuous functions defined on~$\mathfrak {MC}$, associates:
		
		$$
		{\cal{E}} (u,v)= \displaystyle \lim_{m \to + \infty} {\cal{E}}_{m }\left (u_{\mid V_m},v_{\mid V_m}\right)=
		\displaystyle \lim_{m \to + \infty}\displaystyle \sum_{X  \underset{m }{\sim}  Y} r^{-m}\,4^{ 2\,m \,\delta}\,\,\left (u_{\mid V_m}(X)-u_{\mid V_m}(Y)\right )\,\left(v_{\mid V_m}(X)-v_{\mid V_m}(Y)\right) $$

		\noindent satisfies the self-similarity relation:
		
		$$
		{\cal{E}} (u,v)= r ^{-1}\,4^{ 2\,\delta}\,\displaystyle \sum_{i=1}^{8}   {\cal{E}} \left ( u \circ f_i,v\circ f_i \right)$$

	\end{pte}
	
	\vskip 1cm
	
	\begin{proof}

		$$\begin{array}{ccc}  \displaystyle \sum_{i=1}^{8}   {\cal{E}} \left ( u \circ f_i,v\circ f_i \right)&=& \displaystyle \lim_{m \to + \infty} \displaystyle \sum_{i=1}^{8}  {\cal{E}}_{m  }\left ( u_{\mid V_{m }} \circ f_i,v_{\mid V_{m }}\circ f_i \right)\\
		&=&\displaystyle \lim_{m \to + \infty}\displaystyle \sum_{X  \underset{m  }{\sim}  Y}  r^{-m}\,4^{ 2\,m \,\delta}\,  \displaystyle \sum_{i=1}^{8}   \left (u_{\mid V_{m }}\left (f_i(X)\right )-
		u_{\mid V_{m }} \left (f_i(Y)\right  )\right )\,\left(v_{\mid V_{m }}\left (f_i(X)\right )-v_{\mid V_{m }}\left (f_i(Y)\right )\right) \\
		&=&\displaystyle \lim_{m \to + \infty}\displaystyle \sum_{X  \underset{m+1  }{\sim}  Y}  r^{-m}\,4^{ 2\,m \,\delta}\,  \displaystyle \sum_{i=1}^{8}   \left (u_{\mid V_{m+1}} (X) -u_{\mid V_{m+1}}  (Y) \right )\,
		\left(v_{\mid V_{m+1}} (X) -v_{\mid V_{m+1}} (Y) \right) \\
		&=&\displaystyle \lim_{m \to + \infty}  r\,4^{- 2\,\delta}\, {\cal{E}}_{m+1  }\left (u_{\mid V_{m+1}},v_{\mid V_{m+1}} \right) \\
		&=& r\,4^{- 2\,\delta}\, {\cal{E}} (u,v) \\
		\end{array} $$

	\end{proof}
	\vskip 1cm
	
	\begin{notation}
		\noindent We will denote by~$\text{dom}\,{\cal E}$ the subspace of continuous functions defined on~$\mathfrak {MC}$, such that:
		
		$$\mathcal{E}(u)< + \infty$$

	\end{notation}

	\vskip 1cm
	
	\begin{notation}
		\noindent We will denote by~$\text{dom}_0\,{ \cal E}$ the subspace of continuous functions defined on~$\mathfrak {MC}$, which take the value $0$ on~$V_0$, such that:
		
		$$\mathcal{E}(u)< + \infty$$

	\end{notation}

	\vskip 1cm
	\begin{proposition}
		The space $\text{dom} \mathcal{E}$, modulo the space of constant function on $\mathfrak{MC}$, is a Hilbert space.\\
	\end{proposition}
	
	\vskip 1cm
	\newpage
	
	\section{Laplacian, on the Minkowski Curve}
	
	\begin{definition}\textbf{Self-similar measure, on the domain delimited by the Minkowski Curve}\\
		
		\noindent A measure~$\mu$ on~$\R^2$ will be said to be \textbf{self-similar} domain delimited by the Minkowski Curve, if there exists a family of strictly positive pounds~\mbox{$\left (\mu_i\right)_{1 \leq i \leq 8}$} such that:
		
		$$ \mu= \displaystyle \sum_{i=1}^{8} \mu_i\,\mu\circ f_i^{-1} \quad, \quad \displaystyle \sum_{i=1}^{8} \mu_i =1$$
		
		\noindent For further precisions on self-similar measures, we refer to the works of~J.~E.~Hutchinson~(see \cite{Hutchinson1981}).
		
	\end{definition}
	
	\vskip 1cm
	
	\begin{pte}\textbf{Building of a self-similar measure, for the Minkowski Curve}\\
		
		\noindent The Dirichlet forms mentioned in the above require a positive Radon measure with full support.
		
		\noindent Let us set for any integer~$i$ belonging to~$\left \lbrace 1, \hdots, 8 \right \rbrace$:
		
		$$\mu_i=R^{D_{H}\left(\mathfrak {MC}\right)}=\displaystyle \frac{1}{8}$$
		
		\noindent  This enables one to define a self-similar measure~$\mu$ on~$\mathfrak{MC}$ as:
		\[\mu =\displaystyle \frac{1}{8}\,\sum_{i=1}^8  \mu \circ f_i \]

	\end{pte}

	\vskip 1cm
	
	\begin{definition}\textbf{Laplacian of order~$m\,\in\,\N^\star$}\\
		
		\noindent For any strictly positive integer~$m$, and any real-valued function~$u$, defined on the set~$V_m$ of the vertices of the graph~${\mathfrak {MC}}_m $, we introduce the Laplacian of order~$m$,~$\Delta_m(u)$, by:
		
		$$\Delta_m u(X) = \displaystyle\sum_{Y \in V_m,\,Y\underset{m}{\sim} X} \left (u(Y)-u(X)\right)  \quad \forall\, X\,\in\, V_m\setminus V_0 $$

	\end{definition}

	\vskip 1cm

	\begin{definition}\textbf{Harmonic function of order~$m\,\in\,\N^\star$}\\
		
		\noindent Let~$m$ be a strictly positive integer. A real-valued function~$u$,defined on the set~$V_m$ of the vertices of the graph~${\mathfrak {MC}}_m $, will be said to be \textbf{harmonic of order~$m$} if its Laplacian of order~$m$ is null:
		
		$$\Delta_m u(X) =0 \quad \forall\, X\,\in\, V_m\setminus V_0 $$
		
	\end{definition}

	\vskip 1cm
	
	\begin{definition}\textbf{Piecewise harmonic function of order~$m\,\in\,\N^\star$}\\
		
		\noindent  Given a strictly positive integer~$m$, a real valued function~$u$, defined on the set of vertices of~$\mathfrak {MC}$, is said to be \textbf{piecewise harmonic function of order~$m$} if, for any word~${\cal W}$ of length~$ m$,~$u\circ f_{\cal W}$ is harmonic of order~$m$.
		
	\end{definition}
	
	\vskip 1cm

	\begin{definition}\textbf{Existence domain of the Laplacian, for a continuous function on~$\mathfrak {MC}$} (see \cite{Beurling1985})\\
		
		\label{Lapl}
		\noindent We will denote by~$\text{dom}\, \Delta$ the existence domain of the Laplacian, on the graph~${\mathfrak {MC}}$, as the set of functions~$u$ of~$\text{dom}\, \mathcal{E}$ such that there exists a continuous function on~$\mathfrak {MC}$, denoted~$\Delta \,u$, that we will call \textbf{Laplacian of~$u$}, such that :
		$$\mathcal{E}(u,v)=-\displaystyle \int_{{\cal D} \left ( {\mathfrak {MC} }\right)} v\, \Delta u   \,d\mu \quad \text{for any } v \,\in \,\text{dom}_0\, \mathcal{E}$$
	\end{definition}

	\vskip 1cm

	\begin{definition}\textbf{Harmonic function}\\
		
		\noindent A function~$u$ belonging to~\mbox{$\text{dom}\,\Delta$} will be said to be \textbf{harmonic} if its Laplacian is equal to zero.
	\end{definition}
	
	\vskip 1cm

	\begin{notation}
		
		In the following, we will denote by~${\cal H}_0\subset \text{dom}\, \Delta$ the space of harmonic functions, i.e. the space of functions~$u \,\in\,\ \text{dom}\, \Delta$ such that:
		
		$$\Delta\,u=0$$
		
		\noindent Given a natural integer~$m$, we will denote by~${\cal S} \left ({\cal H}_0,V_m \right)$ the space, of dimension~$8^m$, of spline functions " of level~$m$", ~$u$, defined on~${\mathfrak {MC}}$, continuous, such that, for any word~$\cal W$ of length~$m$,~\mbox{$u \circ T_{\cal W}$} is harmonic, i.e.:
		
		$$\Delta_m \, \left ( u \circ T_{\cal W} \right)=0$$
		
	\end{notation}

	\vskip 1cm
	
	\begin{pte}
		
		For any natural integer~$m$:
		
		$${\cal S} \left ({\cal H}_0,V_m \right )\subset  \text{dom }{ \cal E}$$
		
	\end{pte}
	\vskip 1cm

	\begin{pte}
		Let~$m$ be a strictly positive integer,~$X \,\notin\,V_0$ a vertex of the graph~$\mathfrak {MC}$, and~\mbox{$\psi_X^{m}\,\in\,{\cal S} \left ({\cal H}_0,V_m \right)$} a spline  function such that:
		
		$$\psi_X^{m}(Y)=\left \lbrace \begin{array}{ccc}\delta_{XY} & \forall& Y\,\in \,V_m \\
		0 & \forall& Y\,\notin \,V_m \end{array} \right. \quad,  \quad \text{where} \quad    \delta_{XY} =\left \lbrace \begin{array}{ccc}1& \text{if} & X=Y\\ 0& \text{else} &  \end{array} \right.$$

		\noindent Then, since~$X\, \notin \,V_0$: $\psi_X^{m}\,\in \,\text{dom}_0\, \mathcal{E}$.

		\noindent For any function~$u$ of~$\text{dom}\, \mathcal{E}$, such that its Laplacian exists, definition (\ref{Lapl}) applied to~$\psi_X^{m}$ leads to:

		$$\mathcal{E}(u,\psi_X^{m})=\mathcal{E}_m(u,\psi_X^{m})= -r^{-m}\,\Delta_m u(X)=- \displaystyle\int_{{\cal D}\left ({\mathfrak {MC}}\right)}  \psi_X^{m}\,\Delta u  \, d\mu \approx -\Delta  u(X)\, \displaystyle\int_{{\cal D} \left ( {\mathfrak {MC}} \right)}  \psi_X^{m}\, d\mu$$
		
		\noindent since~$\Delta u$ is continuous on~$ {\mathfrak {MC}}$, and the support of the spline function~$\psi_X^{m}$ is close to~$X$:
		
		$$\displaystyle\int_{{\cal D} \left ( {\mathfrak {MC}}  \right)}  \psi_X^{m}\,\Delta u  \, d\mu \approx -\Delta  u(X)\, \displaystyle\int_{{\cal D} \left ( {\mathfrak {MC}} \right)}  \psi_X^{m}\, d\mu$$
		
		\noindent By passing through the limit when the integer~$m$ tends towards infinity, one gets:
		
		$$ \displaystyle \lim_{m \to + \infty} \displaystyle\int_{{\cal D} \left ( {\mathfrak {MC}} \right)}  \psi_X^{m}\,\Delta_m u  \, d\mu=
		\Delta  u(X)\,\displaystyle \lim_{m \to + \infty}   \displaystyle\int_{{\cal D} \left ( {\mathfrak {MC}} \right)}  \psi_X^{m}\, d\mu$$
		
		\noindent i.e.:
		
		$$  \Delta  u(X)= \displaystyle \lim_{m \to + \infty} r^{-m}\,4^{ 2\,m \,\delta}\, \left (  \displaystyle\int_{{\cal D} \left ( {\mathfrak {MC}}  \right)}  \psi_X^{m}\, d\mu  \right)^{-1} \,\Delta_m u(X)\,$$

	\end{pte}

	\vskip 1cm
	
	\begin{remark}

		\noindent As it is explained in~\cite{StrichartzLivre2006}, one has just to reason by analogy with the dimension~1, more particulary, the unit interval~$I=[0,1]$, of extremities~$X_0=(0,0)$, and~$X_1=(1,0)$. The functions~$\psi_{X_1}$ and~$\psi_{X_2}$ such that, for any~$Y$ of~$\R^2$ :
		
		$$\psi_{X_1} (Y)=\delta_{X_1Y} \quad  ,  \quad \psi_{X_2} (Y)=\delta_{X_2Y}   $$
		
		\noindent are, in the most simple way, tent functions. For the standard measure, one gets values that do not depend on~$X_1$, or~$X_2$ (one could, also, choose to fix~$X_1$ and~$X_2$ in the interior of~$I$) :
		
		$$\displaystyle\int_{I}  \psi_{X_1}\, d\mu =\displaystyle\int_{I}  \psi_{X_2}\, d\mu=\displaystyle \frac{1}{2}$$
		
		\noindent (which corresponds to the surfaces of the two tent triangles.) \\

		\begin{figure}[h!]
			\center{\psfig{height=8cm,width=10cm,angle=0,file=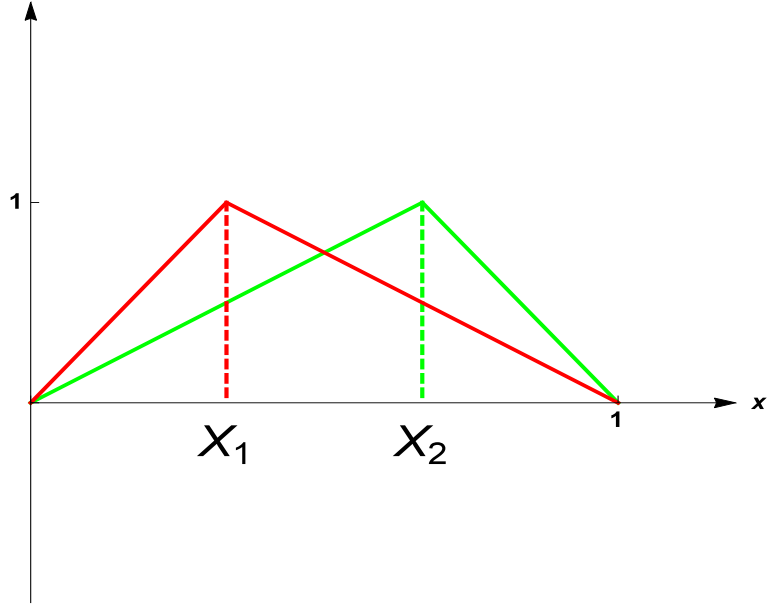}}\\
			\caption{The graphs of the spline functions~$\psi_{X_1}$ and~$\psi_{X_2}$.}
		\end{figure}

		\noindent In our case, we have to build the pendant, we no longer reason on the unit interval, but on our triangular or square cells. \\

		\noindent Given a natural integer~$m$, and a point~$X \,\in\,V_m$, the spline function~$\psi_X^{m}$ is supported by two~$m$-cells. It is such that, for every $m$-line cell $f_{\mathcal{W}} \left( \mathfrak{MC}\right)$ the vertices of which are~$X$,~\mbox{$Y\neq X$} :
		$$\psi_X^{m}+\psi_Y^{m}=1$$

		\noindent Thus:
		\[ \displaystyle \int_{f_{\mathcal{W}} \left( \mathfrak{MC}\right)} \left(\psi_X^{m}+\psi_Y^{m}\right)\, d\mu=\mu(f_{\mathcal{W}} \left( \mathfrak{MC}\right) )=\displaystyle \frac{1}{8^m} \]
		\noindent By symmetry, all three summands have the same integral. This yields:
		
		$$\displaystyle{\int_{f_{\mathcal{W}} \left( \mathfrak{MC}\right)} \psi_X^{m} d\mu=\frac{1}{2\times 8^{m}}}$$
		\noindent Taking into account the contributions of the remaining~$m$-square cells, one has:
		
		$$\displaystyle{\int_{\mathfrak{MC}} \psi_X^{m} d\mu=\frac{1}{ 8^{m}}}$$
		
		\noindent which leads to:
		
		$$\displaystyle \left (\int_{\mathfrak{MC}} \psi_X^{m}\, d\mu \right )^{-1}=\displaystyle 8^{m} $$
		
		\noindent Since:
		
		$$ \displaystyle{r^{-m}= 8^m}$$

		\noindent this enables us to obtain the point-wise formula, for~$u\,\in\,\text{dom}\,\Delta$:
		
		$$\forall \,X \,\in\,{\mathfrak {MC}} \,: \quad  \Delta_{\mu} u(X)= \displaystyle\lim_{m\rightarrow +\infty} 64^m\,
		\Delta_m u(X) $$

	\end{remark}
	
	\vskip 1cm
	
	\vskip 1cm
	\begin{theorem}
		
		\noindent Let~$u$ be in~\mbox{$\text{dom}\,\Delta$}. Then, the sequence of functions~$\left (f_m \right)_{m\in\N}$ such that, for any natural integer~$m$, and any~$X$ of~\mbox{$V_\star\setminus V_0$} :
		
		$$f_m(X)=r^{-m}\,4^{ 2\,m \,\delta}\,\left (\displaystyle \int_{{\cal D} \left ({\mathfrak {MC}}\right) }   \psi_{X}^{m}\,d\mu\right)^{-1}\,\Delta_m \,u(X) $$
		
		\noindent  converges uniformly towards~$\Delta\,u$, and, reciprocally, if the sequence of functions~$\left (f_m \right)_{m\in\N}$ converges uniformly towards a continuous function on~\mbox{$V_\star\setminus V_0$}, then:
		
		$$u \,\in\, \text{dom}\,\Delta$$
	\end{theorem}

	\vskip 1cm
	\begin{proof}

		\noindent Let~$u$ be in~\mbox{$\text{dom}\,\Delta$}. Then:

		$$ r^{-m}\,4^{ 2\,m \,\delta}\,\,\left (\displaystyle \int_{{\cal D} \left ({\mathfrak {MC}}\right) }  \psi_{X}^{m}\,d\mu\right)^{-1}\,\Delta_m \,u(X)=
		\displaystyle \frac{\displaystyle \int_{{\cal D} \left ({\mathfrak {MC}}\right) }  \Delta\,u \,\psi_{X}^{m}\,d\mu} {\displaystyle \int_{{\cal D} \left ({\mathfrak {MC}}\right) }   \psi_{X}^{m}\,d\mu}   $$
		
		\noindent Since~$u$ belongs to~\mbox{$\text{dom}\,\Delta$}, its Laplacian~$\Delta\,u$ exists, and is continuous on the graph~${\mathfrak {MC}}$. The uniform convergence of the sequence~$\left (f_m \right)_{m\in\N}$ follows.\\
		
		\noindent Reciprocally, if the sequence of functions~$\left (f_m \right)_{m\in\N}$ converges uniformly towards a continuous function on~\mbox{$V_\star\setminus V_0$}, the, for any natural integer~$m$, and any~$v$ belonging to~\mbox{$\text{dom}_0\,{\cal E}$}:

		$$\begin{array}{ccc} {\cal{E}}_{m }(u,v)
		&=&  \displaystyle \sum_{(X,Y) \,\in \, V_m^2,\,X  \underset{m }{\sim}  Y} r^{-m}\,4^{ 2\,m \,\delta}\,\,\left (u_{\mid V_m}(X)-u_{\mid V_m}(Y)\right )\,\left(v_{\mid V_m}(X)-v_{\mid V_m}(Y)\right) \\
		&=& \displaystyle \sum_{(X,Y) \,\in \, V_m^2,\,X  \underset{m }{\sim}  Y} r^{-m}\,4^{ 2\,m \,\delta}\,\,\left (u_{\mid V_m}(Y)-u_{\mid V_m}(X )\right )\,\left(v_{\mid V_m}(Y)-v_{\mid V_m}(X)\right) \\
		&=&- \displaystyle \sum_{X \,\in \,V_m\setminus V_0 } r^{-m}\,4^{ 2\,m \,\delta}\,\,\sum_{Y\,\in \,V_m, \, Y  \underset{m }{\sim}  X} v_{\mid V_m}(X)\,\left (u_{\mid V_m}(Y)-u_{\mid V_m}(X)\right )   \\
		& &- \displaystyle \sum_{X \,\in \,  V_0 } r^{-m}\,4^{ 2\,m \,\delta}\,  \sum_{Y\,\in \,V_m ,\, Y  \underset{m }{\sim}  X} v_{\mid V_m}(X)\,\left (u_{\mid V_m}(Y)-u_{\mid V_m}(X)\right )   \\
		&=&- \displaystyle \sum_{X \,\in \,V_m\setminus V_0 } r^{-m}\,4^{ 2\,m \,\delta}\,\,v(X)\,\Delta_m \,u(X)  \\
		&=&- \displaystyle \sum_{X \,\in \,V_m\setminus V_0 } v(X)\,   \left (\displaystyle \int_{{\cal D} \left ({\mathfrak {MC}}\right) }   \psi_{X}^{m}\,d\mu\right) \, r^{-m}\,4^{ 2\,m \,\delta}\, \left (\displaystyle \int_{{\cal D} \left ({\mathfrak {MC}}\right) }   \psi_{X}^{m}\,d\mu\right)^{-1}\, \Delta_m \,u(X)  \\
		\end{array}
		$$
		
		\noindent Let us note that any~$X$ of~$V_m\setminus V_0$ admits exactly two adjacent vertices which belong to~$V_m\setminus V_0$, which accounts for the fact that the sum
		
		$$\displaystyle \sum_{X \,\in \,V_m\setminus V_0 } r^{-m}\,4^{ 2\,m \,\delta}\,\sum_{Y\,\in \,V_m\setminus V_0 ,\, Y  \underset{m }{\sim}  X} v(X)\,\left (u_{\mid V_m}(Y)-u_{\mid V_m}(X)\right)$$
		\noindent has the same number of terms as:
		
		$$ \displaystyle \sum_{(X,Y) \,\in \,(V_m\setminus V_0)^2,\,X  \underset{m }{\sim}  Y} r^{-m}\,8^{ 2\,m \,\delta}\,\left (u_{\mid V_m}(Y)-u_{\mid V_m}(X)\right )\,\left(v_{\mid V_m}(Y)-v_{\mid V_m}(X)\right)  $$
		
		\noindent For any natural integer~$m$, we introduce the sequence of functions~$\left (f_m \right)_{m\in\N}$ such that, for any~$X$ of~$V_m\setminus V_0$:
		
		$$f_m(X)=r^{-m}\,4^{ 2\,m \,\delta}\,\left (\displaystyle \int_{{\cal D} \left ({\mathfrak {MC}} \right) }  \psi_{X}^{m}\,d\mu\right)^{-1}\,\Delta_m \,u(X) $$
		
		\noindent The sequence~$\left (f_m \right)_{m\in\N}$ converges uniformly towards~$\Delta\,u$. Thus:

		$$\begin{array}{ccc} {\cal{E}}_{m }(u,v)
		&=&      -   \displaystyle \int_{{\cal D} \left ({\mathfrak {MC}}\right) }  \left \lbrace \displaystyle \sum_{X \,\in \,V_m\setminus V_0 } v_{\mid V_m}(X)\,\Delta\,u_{\mid V_m}(X)\, \psi_{X}^{m}\right \rbrace \,d\mu
		\end{array}
		$$
		
	\end{proof}

	\section{Normal derivatives}
	
	Let us go back to the case of a function~$u$ twice differentiable on~$I=[0,1]$, that does not vanish in~0 and~:
	$$  \displaystyle \int_0^1 \left (\Delta u \right)(x)\,v(x)\,dx=-   \displaystyle \int_0^1 u'(x)\,v'(x)\,dx+ u '(1)\,v (1)-u '(0)\,v (0) $$
	
	\noindent The normal derivatives:
	
	$$\partial_n u(1)=u'(1) \quad \text{et} \quad \partial_n u(0)=u'(0) $$
	
	\noindent appear in a natural way. This leads to:

	$$  \displaystyle \int_0^1 \left (\Delta u \right)(x)\,v(x)\,dx=-   \displaystyle \int_0^1 u'(x)\,v'(x)\,dx+ \displaystyle \sum_{\partial\, [0,1]} v\,\partial_n\,u $$

	\noindent One meets thus a particular case of the Gauss-Green formula, for an open set~$\Omega$ of~$\R^d$,~$d \,\in\,\N^\star$:

	$$  \displaystyle \int_\Omega \nabla\,  u\, \nabla \, v \, d \mu= -\displaystyle \int_\Omega  \left (\Delta u \right) \,v \,d\mu + \displaystyle \int_{\partial\, \Omega } v\,\partial_n\,u \,d\sigma$$
	
	\noindent where~$\mu$ is a measure on~$\Omega $, and where~$d\sigma$ denotes the elementary surface on~$\partial\, \Omega $.\\
	
	\noindent In order to obtain an equivalent formulation in the case of the graph~$\mathfrak {MC} $, one should have, for a pair of functions~$(u,v)$ continuous on~$\mathfrak {MC} $ such that~$u$ has a normal derivative:

	$$  {\cal E}(u,v)= -\displaystyle \int_\Omega  \left (\Delta u \right) \,v \,d\mu + \displaystyle \sum_{V_0} v\,\partial_n\,u  $$
	
	\noindent For any natural integer~$m$ :

	$$\begin{array}{ccc}
	{\cal{E}}_{m }(u,v)\\
	&=&  \displaystyle \sum_{(X,Y) \,\in \, V_m^2,\,X  \underset{m }{\sim}  Y} r^{-m}\,4^{ 2\,m \,\delta}\,\left (u_{\mid V_m}(Y)-u_{\mid V_m}(X)\right )\,\left(v_{\mid V_m}(Y)-v_{\mid V_m}(X)\right) \\
	&=&- \displaystyle \sum_{X \,\in \,V_m\setminus V_0 } r^{-m}\,8^{ 2\,m \,\delta}\,\sum_{Y\,\in \,V_m  ,\, Y  \underset{m }{\sim}  X} v_{\mid V_m}(X)\,\left (u_{\mid V_m}(Y)-u_{\mid V_m}(X)\right )   \\
	&&-  \displaystyle \sum_{X \,\in \,  V_0 } r^{-m}\,4^{ 2\,m \,\delta}\,\sum_{Y\,\in \,  V_m ,\, Y  \underset{m }{\sim}  X} v_{\mid V_m}(X)\,\left (u_{\mid V_m}(Y)-u_{\mid V_m}(X)\right )  \\
	&=&- \displaystyle \sum_{X \,\in \,V_m\setminus V_0 } v_{\mid V_m}(X)\,r^{-m}\,4^{ 2\,m \,\delta}\,\Delta_m \,u_{\mid V_m}(X)  \\
	&&+ \displaystyle \sum_{X \,\in \,  V_0 } \sum_{Y\,\in \,  V_m ,\, Y  \underset{m }{\sim}  X} r^{-m}\,8^{ 2\,m \,\delta}\,v_{\mid V_m}(X)\,\left (u_{\mid V_m}(X)-u_{\mid V_m}(Y)\right )  \\
	\end{array}
	$$
	
	\noindent We thus come across an analogous formula of the Gauss-Green one, where the role of the normal derivative is played by:

	$$ \displaystyle \sum_{X \,\in \,  V_0 } r^{-m}\,4^{ 2\,m \,\delta}\,\sum_{Y\,\in \,  V_m ,\, Y  \underset{m }{\sim}  X} \,\left (u_{\mid V_m}(X)-u_{\mid V_m}(Y)\right ) $$
	
	\vskip 1cm
	
	\begin{definition}
		
		\noindent For any~$X$ of~$V_0$, and any continuous function~$u$ on~$\mathfrak {MC} $, we will say that~$u$ admits a normal derivative in~$X$, denoted by~$\partial_n\,u(X)$, if:
		
		$$ \displaystyle \lim_{m \to + \infty}  r^{-m}\,4^{ 2\,m \,\delta}\,\sum_{Y\,\in \,  V_m ,\, Y  \underset{m }{\sim}  X} \,\left (u_{\mid V_m}(X)-u_{\mid V_m}(Y)\right ) < + \infty $$
		
		\noindent We will set:
		
		$$\partial_n\,u(X) = \displaystyle \lim_{m \to + \infty}  r^{-m}\,4^{ 2\,m \,\delta}\,\sum_{Y\,\in \,  V_m ,\, Y  \underset{m }{\sim}  X} \,\left (u_{\mid V_m}(X)-u_{\mid V_m}(Y)\right ) < + \infty $$
		
	\end{definition}

	\vskip 1cm
	
	\begin{definition}
		
		\noindent For any natural integer~$m$, any~$X$ of~$V_m$, and any continuous function~$u$ on~$\mathfrak {MC} $, we will say that~$u$ admits a normal derivative in~$X$, denoted by~$\partial_n\,u(X)$, if:
		$$ \displaystyle \lim_{k \to + \infty}  r^{-k}\,4^{ 2\,k \,\delta}\,\sum_{Y\,\in \,  V_k ,\, Y  \underset{k }{\sim}  X} \,\left (u_{\mid V_k} (X)-u_{\mid V_k}(Y)\right ) < + \infty $$
		
		\noindent We will set:
		
		$$\partial_n\,u(X) = \displaystyle \lim_{k \to + \infty} r^{-k}\,4^{ 2\,k \,\delta}\, \sum_{Y\,\in \,  V_k ,\, Y  \underset{k }{\sim}  X} \,\left (u_{\mid V_k}(X)-u_{\mid V_k}(Y)\right ) < + \infty $$
		
	\end{definition}

	\vskip 1cm
	\begin{remark} One can thus extend the definition of the normal derivative of~$u$ to~$\mathfrak {MC} $.
	\end{remark}

	\vskip 1cm
	
	\begin{theorem}

		\noindent Let~$u$ be in~\mbox{$\text{dom}\,\Delta$}. The, for any~$X$ of~$\mathfrak {MC} $,~$\partial_n\,u(X)$ exists. Moreover, for any~$v$ of~\mbox{$\text{dom}\,\cal E$}, et any natural integer~$m$, the Gauss-Green formula writes:

		$$  {\cal E}(u,v) = -\displaystyle \int_{{\cal D} \left ({\mathfrak {MC}}\right ) } \left (\Delta u \right) \,v \,d\mu + \displaystyle \sum_{V_0} v\,\partial_n\,u  $$
	\end{theorem}
	\vskip 1cm

	\section{Spectrum of the Laplacian}

	\subsection{Spectral decimation}

	In the following, let~$u$ be in~$\text{dom}\, \Delta$. We will apply the \emph{\textbf{spectral decimation method}} developed by~R.~S.~Strichartz \cite{StrichartzLivre2006}, in the spirit of the works of M.~Fukushima et T.~Shima \cite{Fukushima1994}. In order to determine the eigenvalues of the Laplacian~$\Delta\, u$ built in the above, we concentrate first on the eigenvalues~$\left (-{\lambda_m}\right)_{m\in\N}$ of the sequence of graph Laplacians~$\left (\Delta_m \,u\right)_{m\in\N}$, built on the discrete sequence of graphs~$\left ({\mathfrak {MC}}_m\right)_{m\in\N}$. For any natural integer~$m$, the restrictions of the eigenfunctions of the continuous Laplacian~$\Delta\,u$ to the graph~${\mathfrak {MC}}_m$ are, also, eigenfunctions of the Laplacian~$\Delta_m$, which leads to recurrence relations between the eigenvalues of order~$m$ and~$m+1$.
	
	\vskip 1cm
	
	We thus aim at determining the solutions of the eigenvalue equation:

	$$-\Delta\,u=\lambda\,u \quad  \text{on }{ \mathfrak {MC}}$$
	
	\noindent as limits, when the integer~$m$ tends towards infinity, of the solutions of:

	$$-\Delta_m\,u=\lambda_m\,u \quad  \text{on }V_m\setminus V_0$$
	
	\noindent We will call them Dirichlet eigenvalues (resp. Neumann eigenvalues) if:
	\[u_{\mid \partial{ \mathfrak {MC}}}=0 \quad \left( \text{resp.} \quad {\partial_n u}_{\mid \partial{ \mathfrak {MC}}}=0\right) \]

	\vskip 1cm
	\noindent Given a strictly positive integer~$m$, let us consider a~$(m-1)-$cell, with boundary vertices~$X_0$,~$X_1$.\\
	\noindent We denote by~$Y_1$,~$Y_2$,~$Y_3$,~$Y_4$,~$Y_5$,~$Y_6$,~$Y_7$ the points of $V_m\setminus V_{m-1}$ (see Fig. ):
	
	\begin{center}
		\includegraphics[scale=1]{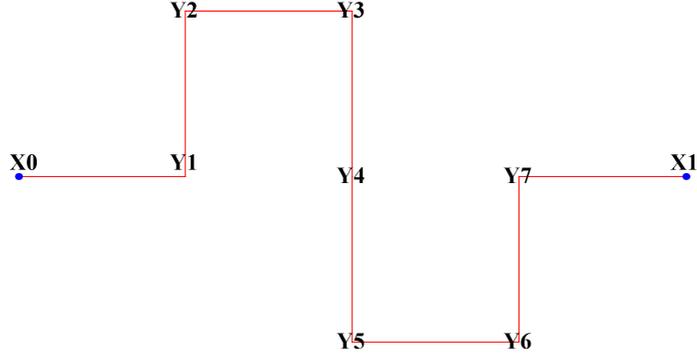}
		\captionof{figure}{The $m$-cell}
		\label{fig10}
	\end{center}

	\noindent  The discrete equation on $\mathfrak{MC}$ leads to the following system:
	
	$$\left \lbrace  \begin{array}{ccc}
	\left(2-\lambda_m \right)u(Y_1)&=&u(X_0)+u(Y_2)\\
	\left(2-\lambda_m \right)u(Y_2)&=&u(Y_1)+u(Y_3)\\
	\left(2-\lambda_m \right)u(Y_3)&=&u(Y_2)+u(Y_4)\\
	\left(2-\lambda_m \right)u(Y_4)&=&u(Y_3)+u(Y_5)\\
	\left(2-\lambda_m \right)u(Y_5)&=&u(Y_4)+u(y_6)\\
	\left(2-\lambda_m \right)u(Y_6)&=&u(Y_5)+u(Y_7)\\
	\left(2-\lambda_m \right)u(Y_7)&=&u(Y_6)+u(X_1)
	\end{array}  \right.$$
	
	Under matrix form, one has:
	
	$$  \mathbf{A}_m \mathbf{u}= \mathbf{u}_0$$
	
	where:
	
	$$ \mathbf{A}_m = \left( \begin{matrix}
	2-\lambda_m & -1 & 0 & 0 & 0 & 0 & 0 \\
	-1 & 2-\lambda_m & -1 & 0 & 0 & 0 & 0 \\
	0 & -1 & 2-\lambda_m & -1 & 0 & 0 & 0 \\
	0 & 0 & -1 & 2-\lambda_m & -1 & 0 & 0 \\
	0 & 0 & 0 & -1 & 2-\lambda_m & -1 & 0 \\
	0 & 0 & 0 & 0 & -1 & 2-\lambda_m & -1 \\
	0 & 0 & 0 & 0 & 0 & -1 & 2-\lambda_m \\
	\end{matrix}\right) $$
	
	$$\mathbf{u}=\left( \begin{matrix}
	u(Y_1)\\
	u(Y_2)\\
	u(Y_3)\\
	u(Y_4)\\
	u(Y_5)\\
	u(Y_6)\\
	u(Y_7)\\
	\end{matrix} \right) \ quad , \quad \mathbf{u}_0=\left( \begin{matrix}
	u(X_0)\\
	0\\
	0\\
	0\\
	0\\
	0\\
	u(X_1)\\
	\end{matrix} \right) $$
	
	\newpage
	\subsubsection{Direct method}

	\noindent By assuming~$\lambda_m \neq \left\{2,2+\varepsilon\, \sqrt{2}, 2+\varepsilon\,\sqrt{2+\varepsilon\,\sqrt{2}} \right\}$,~\mbox{$\varepsilon\,\in\,\left \lbrace -1,1 \right \rbrace $} for which the matrix $\mathbf{A}_m$ is not invertible, we can solve the system using Gauss algorithm applied on the augmented matrix $\tilde{\mathbf{A}}_m$ :
  	\begin{align*}
 	 	\tilde{\mathbf{A}}_m = & \left( \begin{matrix}
 	 	2-\lambda_m & -1 & 0 & 0 & 0 & 0 & 0 & u(X_0)\\
 	 	-1 & 2-\lambda_m & -1 & 0 & 0 & 0 & 0 & 0\\
 	 	0 & -1 & 2-\lambda_m & -1 & 0 & 0 & 0 & 0\\
 	 	0 & 0 & -1 & 2-\lambda_m & -1 & 0 & 0 & 0\\
 	 	0 & 0 & 0 & -1 & 2-\lambda_m & -1 & 0 & 0\\
 	 	0 & 0 & 0 & 0 & -1 & 2-\lambda_m & -1 & 0\\
 	 	0 & 0 & 0 & 0 & 0 & -1 & 2-\lambda_m & u(X_1)\\
 	 	\end{matrix}\right) .
 	 		\end{align*}

	\normalsize
	
Thus:
	
	\begin{align*}
	\left(
	\begin{matrix}
	u(Y_1)\\
	u(Y_2)\\
	u(Y_3)\\
	u(Y_4)\\
	u(Y_5)\\
	u(Y_6)\\
	u(Y_7)\\
	\end{matrix}
	\right)=
	\left(
	\begin{array}{c}
	-\frac{u(X_0) \left(\lambda _m-4\right) \lambda _m \left(\left(\lambda _m-4\right)
		\lambda _m \left(\left(\lambda _m-4\right) \lambda _m+7\right)+14\right)+7
		u(X_0)+u(X_1)}{\left(\lambda _m-2\right) \left(\left(\lambda _m-4\right) \lambda
		_m+2\right) \left(\left(\lambda _m-4\right) \lambda _m \left(\lambda
		_m-2\right){}^2+2\right)} \\
	\frac{u(X_0) \left(\lambda _m-4\right) \lambda _m \left(\lambda _m-2\right){}^2+3
		u(X_0)+u(X_1)}{\left(\left(\lambda _m-4\right) \lambda _m+2\right) \left(\left(\lambda
		_m-4\right) \lambda _m \left(\lambda _m-2\right){}^2+2\right)} \\
	-\frac{\left(\lambda _m-4\right) \lambda _m \left(a \left(\lambda _m-4\right)
		\lambda _m+5 u(X_0)+u(X_1)\right)+5 u(X_0)+3 u(X_1)}{\left(\lambda _m-2\right) \left(\left(\lambda
		_m-4\right) \lambda _m+2\right) \left(\left(\lambda _m-4\right) \lambda _m
		\left(\lambda _m-2\right){}^2+2\right)} \\
	\frac{u(X_0)+u(X_1)}{\left(\lambda _m-4\right) \lambda _m \left(\lambda _m-2\right){}^2+2} \\
	-\frac{\left(\lambda _m-4\right) \lambda _m \left(u(X_0)+u(X_1) \left(\lambda _m-4\right)
		\lambda _m+5 u(X_1)\right)+3 u(X_0)+5 u(X_1)}{\left(\lambda _m-2\right) \left(\left(\lambda
		_m-4\right) \lambda _m+2\right) \left(\left(\lambda _m-4\right) \lambda _m
		\left(\lambda _m-2\right){}^2+2\right)} \\
	\frac{u(X_0)+u(X_1) \left(\lambda _m-4\right) \lambda _m \left(\lambda _m-2\right){}^2+3
		u(X_1)}{\left(\left(\lambda _m-4\right) \lambda _m+2\right) \left(\left(\lambda
		_m-4\right) \lambda _m \left(\lambda _m-2\right){}^2+2\right)} \\
	-\frac{u(X_0)+u(X_1) \left(\lambda _m-4\right) \lambda _m \left(\left(\lambda _m-4\right)
		\lambda _m \left(\left(\lambda _m-4\right) \lambda _m+7\right)+14\right)+7
		u(X_1)}{\left(\lambda _m-2\right) \left(\left(\lambda _m-4\right) \lambda _m+2\right)
		\left(\left(\lambda _m-4\right) \lambda _m \left(\lambda
		_m-2\right){}^2+2\right)} \\
	\end{array}
	\right).
	\end{align*}

	\noindent Let us now compare the $\lambda_{m-1}-$eigenvalues on~$V_{m-1}$, and the $\lambda_{m}-$eigenvalues on~$V_{m}$. To this purpose, we fix~\mbox{$X_0 \,\in V_{m-1} \setminus V_0$}.\\
	\noindent One has to bear in mind that~$X_0$ also belongs to a~$(m-1)-$cell, with boundary points $X_0$,~$X'_1$ and interior points~$Z_1$,~$Z_2$,~$Z_3$,~$Z_4$,~$Z_5$,~$Z_6$,~$Z_7$.\\
	\noindent Thus:
	
	$$ \left(2-\lambda_{m-1}
	\right)u(X_0)=u(X_1)+u(X'_1)  \quad ,\quad
	\left(2-\lambda_{m}
	\right)u(X_0) =u(Y_1)+u(Z_1)
	$$
	
	\noindent and:
	\begin{align*}
	u(Y_1)+u(Z_1)=\left\{-\frac{2 u(X_0) ((\lambda_m-3) (\lambda_m-2) \lambda_m-1) (\lambda_m ((\lambda_m-7) \lambda_m+14)-7)+u(X_1)+u(X'_1)}{(\lambda_m-2) ((\lambda_m-4) \lambda_m+2)
		\left((\lambda_m-4) \lambda_m (\lambda_m-2)^2+2\right)}\right\}
	\\
	\end{align*}

	\noindent Then :
	
	\begin{align*}
	\left(2-\lambda_{m}
	\right)u(X_0)&=\left\{-\frac{2 u(X_0) ((\lambda_m-3) (\lambda_m-2) \lambda_m-1) (\lambda_m ((\lambda_m-7) \lambda_m+14)-7)+u(X_1)+u(X'_1)}{(\lambda_m-2) ((\lambda_m-4) \lambda_m+2)
		\left((\lambda_m-4) \lambda_m (\lambda_m-2)^2+2\right)}\right\}
	\\
	&=u(X_0)\left\{-\frac{2  ((\lambda_m-3) (\lambda_m-2) \lambda_m-1) (\lambda_m ((\lambda_m-7) \lambda_m+14)-7)+\left(2-\lambda_{m-1}
		\right)}{(\lambda_m-2) ((\lambda_m-4) \lambda_m+2)
		\left((\lambda_m-4) \lambda_m (\lambda_m-2)^2+2\right)}\right\} \\
	\end{align*}

	\noindent Finally:
	
	\[
	\lambda_{m-1}=-(\lambda_m-4) (\lambda_m-2)^2 \lambda_m ((\lambda_m-4) \lambda_m + 2)^2
	\]
	
	\noindent One may solve:
	
	$$
	\lambda_m =     2+\varepsilon_1\, \sqrt{2+\varepsilon_2\, \sqrt{2+\varepsilon_3\, \sqrt{4-\lambda_{m-1}}}}
	\quad, \quad  \varepsilon_1,\varepsilon_2,\varepsilon_2,\,\in\,\left \lbrace -1,1 \right \rbrace
	$$

	\noindent Let us introduce:
	
	\[\lambda = \displaystyle \lim_{m\rightarrow\infty}r^{-m}\,4^{ 2\,m \,\delta}\,\left (\displaystyle \int_{{\cal D} \left ({\mathfrak {MC}} \right) }  \psi_{X}^{m}\,d\mu\right)^{-1}\,\lambda_m
	= \displaystyle \lim_{m\rightarrow\infty} 64^{ m} \,\lambda_m\]
	
	\noindent One may note that the limit exists, since, when~$x$ is close to~0:
	
	$$ 2-\sqrt{2+\sqrt{2+\sqrt{4-x}}}=\frac{x}{64}+O\left(x^2\right)$$

	\subsubsection{Recursive method}
	
	It is interesting to apply the method outlined by the second author ; the sequence~$\left ( u \left ( Y_{k } \right) \right)_{0 \leq i \leq 8}$ satisfies a second order recurrence relation, the characteristic equation of which is:

	
	$$ r^2+\left \lbrace \lambda_{m }-2 \right\rbrace \,r+1=0 $$
	
	\noindent The discriminant is:
	
	$$ \delta_m= \left \lbrace \lambda_{m }-2 \right\rbrace ^2-4= \omega_m^2 \quad, \quad \omega_m\,\in\,\C$$

	\noindent The roots~$r_{1,m}$ and~$r_{2,m}$ of the characteristic equation are the scalar given by:

	$$r_{1,m}=\displaystyle \frac{2-\lambda_{m }-\omega_m}{2}  \quad, \quad r_{2,m}=\displaystyle \frac{2-\lambda_{m }+\omega_m}{2}$$
	
	\noindent One has then, for any natural integer~$k$ of~\mbox{$\left \lbrace 0,\hdots,8 \right \rbrace $} :
	
	$$ u \left ( Y_{k  } \right) = \alpha_m\, r_{1,m}^k +\beta_m\, r_{2,m}^k$$
	
	\noindent where~$\alpha_m$ and~$\beta_m$ denote scalar constants.
	
	\noindent In the same way,~$X_0$ and~$X_1$ belong to a sequence that satisfies a second order recurrence relation, the characteristic equation of which is:

	$$ \left \lbrace \Lambda_{m-1 }-2 \right\rbrace \,r= - 1 -r^2$$
	
	\noindent and discriminant:
	
	$$ \delta_{m-1}= \left \lbrace \Lambda_{m-1 }-2 \right\rbrace ^2-4= \omega_{m-1}^2 \quad, \quad \omega_{m-1}\,\in\,\C$$
	
	\noindent The roots~$r_{1,m-1}$ and~$r_{2,m-1}$ of this characteristic equation are the scalar given by:

	$$r_{1,m-1}=\displaystyle \frac{2-\Lambda_{m-1}-\omega_{m-1}}{2} \quad, \quad r_{2,m-1}=\displaystyle \frac{2-\Lambda_{m-1 }+\omega_{m-1}}{2}$$

	\noindent From this point, the compatibility conditions, imposed by spectral decimation, have to be satisfied:

	$$  \left \lbrace \begin{array}{ccc}
	u \left ( Y_{0  } \right)&=&u \left ( X_{0  } \right)  \\
	u \left ( Y_{8 } \right)&=&u \left ( X_{ 1 } \right)  \\
	\end{array}\right.$$
	
	\noindent i.e.:

	$$  \left \lbrace \begin{array}{ccccc}
	\alpha_{m }  +\beta_{m } &=& \alpha_{m-1}  +\beta_{m-1}  & {\cal C}_{ m}\\
	\alpha_{m }\, r_{1,m }^{8} +\beta_{m }\, r_{2,m }^{8}&=& \alpha_{m-1}\, r_{1,m-1}  +\beta_{m-1}\, r_{2,m-1}  & {\cal C}_{2, m}\\
	\end{array}\right. $$

	\noindent where, for any natural integer~$m$,~$\alpha_{m }$ and~$\beta_{m }$ are scalar constants (real or complex).\\
	
	\noindent Since the graph~$ {\mathfrak MI}_{ m-1} $ is linked to the graph~${\mathfrak MI}_{ m }$ by a similar process to the one that links~${\mathfrak MI}_{  1}$ to~${\mathfrak MI}_{0}$, one can legitimately consider that the constants~$\alpha_m$ and~$\beta_m$ do not depend on the integer~$m$:
	
	$$\forall\,m\,\in\,\N^\star \, : \quad \alpha_m=\alpha \,\in\,\R  \quad, \quad \beta_m=\beta \,\in\,\R $$

	\noindent The above system writes:

	$$
	\alpha \, r_{1,m }^{8} +\beta \, r_{2,m }^{8 } =  \alpha \, r_{1,m-1}  +\beta \, r_{2,m-1}
	$$
	
	\noindent and is satisfied for:

	$$  \left \lbrace \begin{array}{ccc}
	r_{1,m }^{8} &=& r_{1,m -1}     \\
	r_{2,m }^{8}&=&  r_{2,m-1}  \\
	\end{array}\right. $$
	
	\noindent i.e.:

	$$  \left \lbrace \begin{array}{ccc}
	\left ( \displaystyle \frac{2-\lambda_{m }- \omega_m}{2}  \right)^{8} &=&\displaystyle \frac{2-\lambda_{m-1 }-\omega_{m-1}}{2}     \\
	\left ( \displaystyle \frac{2-\lambda_{m }+ \omega_m}{2}  \right)^{8} &=&\displaystyle \frac{2-\lambda_{m-1 }+\omega_{m-1}}{2}     \\
	\end{array}\right. $$
	
	\noindent This lead to the recurrence relation:
	
	$$\forall \,m\,\,\N,\, n \geq 2: \quad
	\lambda_{m-1}=-(\lambda_m -4) (\lambda_m -2)^2 \lambda_m (2 + (\lambda_m-4) \lambda_m)^2
	$$
	
	\noindent which is the same one as in the above.
	

	\subsection{A spectral means of determination of the ramification constant}

	In the sequel, we present an alternative method, that enables one to compute the ramification constant using spectral decimation (we refer to~Zhou~\cite{Zhou2007} for further details). Given a strictly positive integer~$m$, let us denote by~$H_m$ the Laplacian matrix associated to the graph~$\mathfrak{MC}_m$, and by~$\cal L$ the set of linear real valued functions defined on~$\mathfrak{MC}_m$. One may write:

	\[
	H_m=
	\left[
	\begin{matrix}
	T_m & J_m^T \\
	J_m & X_m \\
	\end{matrix}
	\right]
	\]
	
	\noindent where $T_m \,\in \,{\cal L}(V_0)$, $J_m\, \in\ {\cal  L}(V_0,V_m\setminus V_0)$ and $X_m \,\in \,{\cal L}(V_m\setminus V_0)$.\\
	
	Let $D$ be the Laplacian matrix on $\mathfrak{MC}_0$, and $M$ a diagonal matrix with $M_{ii}=-X_{ii}$.\\
	
	\vskip 1cm
	
	\begin{definition}
		The Laplacian is said to have a strong harmonic structure if there exist rational functions of the real variable~$\lambda$, respectively denoted by~$K_D $ and $K_T$ such that, for for any real number~$\lambda$ satisfying~\mbox{$\det \left (X + \lambda\, M\right ) \neq 0$}:
		$$T + J^T \,(X +\lambda M)^{-1}J = K_D(\lambda)\,D + K_T(\lambda)\,T $$

		\noindent Let us set:
		\[\mathbf{F}  =\left \lbrace \lambda \in \R \ : \quad K_D(0) = 0 \, \text{or} \ \det (X +\lambda M) = 0 \right \rbrace  \]
		\noindent and:
		\[\mathbf{F}_k  =\left \lbrace\lambda \in \mathbf{F} \,\big |\,  \lambda \ \text{is an eigenvalue of the normalized laplacian} \right \rbrace \]
		
		\noindent The elements of~$\mathbf{F}$ will be called \textbf{forbidden eigenvalues}. As for the elements of~$\mathbf{F}_k$, they will be called \textbf{forbidden eigenvalues at the~$k^{th}$ step}.
		
		\noindent The map~$R$ such that, for any real number~$\lambda$:
		\[
		R(\lambda)  =\displaystyle \frac{\lambda-K_T(\lambda)}{K_D(\lambda)}
		\]
		\noindent will be called \textbf{spectral decimation function}.
	\end{definition}

	\vskip 1cm
	
	\begin{remark}
		In the case of the Minkowski Curve, one may check that:
		
		\[ \#(F_i(V_0)\cap V_0)\leq 1 \quad \forall \ i \in \{1,2,3,4,5,6,7,8\} \]

		$$
		D=
		\left(
		\begin{array}{cc}
		-1 & 1 \\
		1 & -1 \\
		\end{array}
		\right)
		\quad , \quad
		T=
		\left(
		\begin{array}{cc}
		-1 & 0 \\
		0 & -1 \\
		\end{array}
		\right)
		$$
		
		\[
		X=
		\left(
		\begin{array}{ccccccc}
		-2 & 1 & 0 & 0 & 0 & 0 & 0 \\
		1 & -2 & 1 & 0 & 0 & 0 & 0 \\
		0 & 1 & -2 & 1 & 0 & 0 & 0 \\
		0 & 0 & 1 & -2 & 1 & 0 & 0 \\
		0 & 0 & 0 & 1 & -2 & 1 & 0 \\
		0 & 0 & 0 & 0 & 1 & -2 & 1 \\
		0 & 0 & 0 & 0 & 0 & 1 & -2 \\
		\end{array}
		\right)
		\]
		
		$$
		M=
		\left(
		\begin{array}{ccccccc}
		2 & 0 & 0 & 0 & 0 & 0 & 0 \\
		0 & 2 & 0 & 0 & 0 & 0 & 0 \\
		0 & 0 & 2 & 0 & 0 & 0 & 0 \\
		0 & 0 & 0 & 2 & 0 & 0 & 0 \\
		0 & 0 & 0 & 0 & 2 & 0 & 0 \\
		0 & 0 & 0 & 0 & 0 & 2 & 0 \\
		0 & 0 & 0 & 0 & 0 & 0 & 2 \\
		\end{array}
		\right)
		\quad , \quad
		J=
		\left(
		\begin{array}{cc}
		1 & 0 \\
		0 & 0 \\
		0 & 0 \\
		0 & 0 \\
		0 & 0 \\
		0 & 0 \\
		0 & 1 \\
		\end{array}
		\right)
		$$
		
		\[
		T + J^T \,(X +\lambda M)^{-1}\,J=
		\left(
		\begin{array}{cc}
		-\frac{(2 \lambda -1)\, \left(4 \lambda^2-6 \lambda +1\right) (8 (\lambda -2) (\lambda -1) \lambda  (2 \lambda -3)+1)}{8 (\lambda -1) (2
			(\lambda -2) \lambda +1) \left(8 (\lambda -2) \lambda  (\lambda -1)^2+1\right)} & -\frac{1}{8 (\lambda -1) (2 (\lambda -2) \lambda +1)
			\left(8 (\lambda -2) \lambda  (\lambda -1)^2+1\right)} \\
		-\frac{1}{8 (\lambda -1) (2 (\lambda -2) \lambda +1) \left(8 (\lambda -2) \lambda  (\lambda -1)^2+1\right)} & -\frac{(2 \lambda -1)
			\left(4 \lambda^2-6 \lambda +1\right) (8 (\lambda -2) (\lambda -1) \lambda  (2 \lambda -3)+1)}{8 (\lambda -1) (2 (\lambda -2) \lambda +1)
			\left(8 (\lambda -2) \lambda  (\lambda -1)^2+1\right)} \\
		\end{array}
		\right)
		\]
		\noindent Thus, the Minkowsky Curve has a strong harmonic structure, and:
		\footnotesize
		\[T + J^T \,(X +\lambda  M)^{-1}J=-\frac{1}{8
			(\lambda -1) (2 (\lambda -2) \lambda +1) \left(8 (\lambda -2) \lambda  (\lambda -1)^2+1\right)}D+
		\frac{\lambda  (4 \lambda  (\lambda  (2 \lambda -7)+7)-7)}{8 (\lambda -2) \lambda  (\lambda -1)^2+1}T\]
		\normalsize
		\noindent Then :
		
		\[\frac{1}{K_D(0)}=8\]
		
		\noindent And we can verify again the spectral decimation function is :
		
		\[
		R(\lambda)  =-(\lambda-4) (\lambda-2)^2 \lambda ((\lambda-4) \lambda + 2)^2
		\]

	\end{remark}
	
	\vskip 1cm
	
	\subsection{A detailed study of the spectrum}

	\subsubsection{First case:~$m=1$.}

	\noindent The Minkowski graph, with its eleven vertices, can be seen in the following figures:
	
	\begin{center}
		\includegraphics[scale=1]{MIH2.eps}
		\captionof{figure}{The line after the first iteration.}
		\label{fig1}
	\end{center}

	\noindent Let us look for the kernel of the matrix $\mathbf{A}_1$ in the case of forbidden eigenvalues i.e. $$\lambda_1 \in \left\{2,2+\varepsilon \, \sqrt{2}, 2+\varepsilon\,\sqrt{2+\varepsilon\,\sqrt{2}},
	\right\}
	\quad , \quad \varepsilon \,\in\,\left \lbrace -1, 1 \right \rbrace$$

	\noindent $\rightsquigarrow$ \underline{For~$\lambda_1=2$, we find the one dimensional Dirichlet eigenspace:}
	
	$$V^1_2=\text{Vect}\, \left \lbrace (-1, 0, 1, 0, -1, 0, 1)\right \rbrace $$

	\noindent $\rightsquigarrow$ \underline{For~$\lambda_1=2 + \varepsilon\,\sqrt{2}$,~\mbox{$\varepsilon \,\in\,\left \lbrace -1, 1 \right \rbrace$}, we find the one-dimensional Dirichlet eigenspace:}
	
	$$V^1_{2 + \varepsilon\in\{-1,1\}\sqrt{2}}=\text{Vect}\, \left \lbrace (-1, \varepsilon\,\sqrt{2}, -1, 0, 1, -\varepsilon\,\sqrt{2}, 1
	)\right \rbrace$$

	

	\noindent $\rightsquigarrow$ \underline{For~$\lambda_1=2+\varepsilon \,\in\,\left \lbrace -1, 1 \right \rbrace\sqrt{2+\sqrt{2}}$,~\mbox{$\varepsilon \,\in\,\left \lbrace -1, 1 \right \rbrace$}, we find the one-dimensional Dirichlet eigenspace:}
	
	$$V^1_{2+\varepsilon \,\in\,\left \lbrace -1, 1 \right \rbrace\sqrt{2+\sqrt{2}}}=
	\text{Vect}\, \left \lbrace
	\left(
	\begin{array}{ccccccc}
	1, & -\varepsilon \, \sqrt{2+\sqrt{2}}, & 1+\sqrt{2}, & -\varepsilon \, \sqrt{2 \left(2+\sqrt{2}\right)}, &
	1+\sqrt{2}, & -\varepsilon\,\sqrt{2+\sqrt{2}}, & 1 \\
	\end{array}
	\right)
	\right \rbrace
	$$
	
	
	
	\noindent $\rightsquigarrow$ \underline{ For~$\lambda_1=2+\varepsilon\, \sqrt{2-\sqrt{2}}$,~\mbox{$\varepsilon \,\in\,\left \lbrace -1, 1 \right \rbrace$}, we find the one-dimensional Dirichlet eigenspace}
	
	$$V^1_{2+\varepsilon\in\{-1,1\}\sqrt{2-\sqrt{2}}}=
	\text{Vect}\, \left \lbrace
	\left(
	\begin{array}{ccccccc}
	1, & -\varepsilon\,\sqrt{2-\sqrt{2}}, & 1-\sqrt{2}, & \varepsilon\,\sqrt{2 \left(2-\sqrt{2}\right)}, & 1-\sqrt{2},
	& -\varepsilon\,\sqrt{2-\sqrt{2}}, & 1 \\
	\end{array}
	\right)
	\right \rbrace$$
	
	
	
	\noindent One may easily check that:
	$$\# (V_1\setminus  V_0) =7$$
	
	\noindent Thus, the spectrum is complete.\
	
	\vskip 1cm
	
	\subsubsection{Second case:~$m=2$}
	
	\noindent Let us now move to the~$m=2$ case.
	
	\begin{center}
		\includegraphics[scale=1]{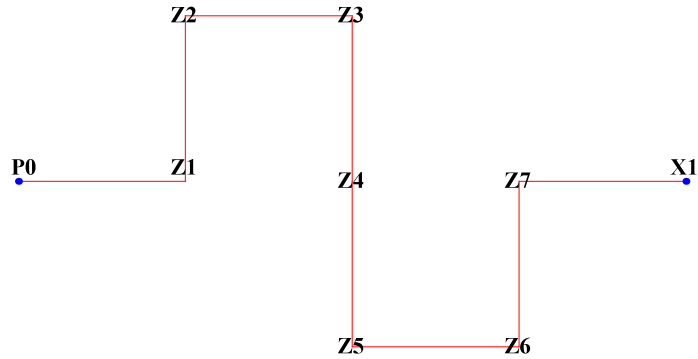}
		\captionof{figure}{The cell $F_1 (V_1)$}
		\label{fig1}
	\end{center}
	
	\noindent Let us denote by $Z^i_j:=f_i(X_j)$ the points of $V_2\setminus V_1$ that belongs to the cell $f_i(V_0)$.\\
	
	\noindent One has to solve the following systems, taking into account the Dirichlet boundary conditions~\mbox{$\left (u(P_0)=u(P_1)=0\right)$}:
	
	$$\left \lbrace \begin{array}{ccc}
	\left(2-\lambda_m \right)u(X_1)&=&u(Z^1_7)+u(Z^2_1)\\
	\left(2-\lambda_m \right)u(X_2)&=&u(Z^2_7)+u(Z^3_1)\\
	\left(2-\lambda_m \right)u(X_3)&=&u(Z^3_7)+u(Z^4_1)\\
	\left(2-\lambda_m \right)u(X_4)&=&u(Z^4_7)+u(Z^5_1)\\
	\left(2-\lambda_m \right)u(X_5)&=&u(Z^5_7)+u(Z^6_1)\\
	\left(2-\lambda_m \right)u(X_6)&=&u(Z^6_7)+u(Z^7_1)\\
	\left(2-\lambda_m \right)u(X_7)&=&u(Z^7_7)+u(Z^8_1)\\
	\end{array}\right. \quad , \quad
	\left \lbrace \begin{array}{ccc}
	\left(2-\lambda_m \right)u(Z^1_1)&=&u(X_0)+u(Z^1_2)\\
	\left(2-\lambda_m \right)u(Z^1_2)&=&u(Z^1_1)+u(Z^1_3)\\
	\left(2-\lambda_m \right)u(Z^1_3)&=&u(Z^1_2)+u(Z^1_4)\\
	\left(2-\lambda_m \right)u(Z^1_4)&=&u(Z^1_3)+u(Z^1_5)\\
	\left(2-\lambda_m \right)u(Z^1_5)&=&u(Z^1_4)+u(Z^1_6)\\
	\left(2-\lambda_m \right)u(Z^1_6)&=&u(Z^1_5)+u(Z^1_7)\\
	\left(2-\lambda_m \right)u(Z^1_6)&=&u(Z^1_6)+u(X_1)\\
	\end{array}\right. $$

	$$\left \lbrace \begin{array}{ccc}
	\left(2-\lambda_m \right)u(Z^2_1)&=&u(X_1)+u(Z^2_2)\\
	\left(2-\lambda_m \right)u(Z^2_2)&=&u(Z^2_1)+u(Z^2_3)\\
	\left(2-\lambda_m \right)u(Z^2_3)&=&u(Z^2_2)+u(Z^2_4)\\
	\left(2-\lambda_m \right)u(Z^2_4)&=&u(Z^2_3)+u(Z^2_5)\\
	\left(2-\lambda_m \right)u(Z^2_5)&=&u(Z^2_4)+u(Z^2_6)\\
	\left(2-\lambda_m \right)u(Z^2_6)&=&u(Z^2_5)+u(Z^2_7)\\
	\left(2-\lambda_m \right)u(Z^2_7)&=&u(Z^2_6)+u(X_2)\\
	\end{array}\right. \quad , \quad
	\left \lbrace \begin{array}{ccc}
	\left(2-\lambda_m \right)u(Z^3_1)&=&u(X_2)+u(Z^3_2)\\
	\left(2-\lambda_m \right)u(Z^3_2)&=&u(Z^3_1)+u(Z^3_3)\\
	\left(2-\lambda_m \right)u(Z^3_3)&=&u(Z^3_2)+u(Z^3_4)\\
	\left(2-\lambda_m \right)u(Z^3_4)&=&u(Z^3_3)+u(Z^3_5)\\
	\left(2-\lambda_m \right)u(Z^3_5)&=&u(Z^3_4)+u(Z^3_6)\\
	\left(2-\lambda_m \right)u(Z^3_6)&=&u(Z^3_5)+u(Z^3_7)\\
	\left(2-\lambda_m \right)u(Z^3_7)&=&u(Z^3_6)+u(X_3)\\
	\end{array}\right. $$

	$$
	\left \lbrace \begin{array}{ccc}
	\left(2-\lambda_m \right)u(Z^4_1)&=&u(X_3)+u(Z^4_2)\\
	\left(2-\lambda_m \right)u(Z^4_2)&=&u(Z^4_1)+u(Z^4_3)\\
	\left(2-\lambda_m \right)u(Z^4_3)&=&u(Z^4_2)+u(Z^4_4)\\
	\left(2-\lambda_m \right)u(Z^4_4)&=&u(Z^4_3)+u(Z^4_5)\\
	\left(2-\lambda_m \right)u(Z^4_5)&=&u(Z^4_4)+u(Z^4_6)\\
	\left(2-\lambda_m \right)u(Z^4_6)&=&u(Z^4_5)+u(Z^4_7)\\
	\left(2-\lambda_m \right)u(Z^4_7)&=&u(Z^4_6)+u(X_4)\\
	\end{array}\right. \quad , \quad
	\left \lbrace \begin{array}{ccc}
	\left(2-\lambda_m \right)u(Z^5_1)&=&u(X_4)+u(Z^5_2)\\
	\left(2-\lambda_m \right)u(Z^5_2)&=&u(Z^5_1)+u(Z^5_3)\\
	\left(2-\lambda_m \right)u(Z^5_3)&=&u(Z^5_2)+u(Z^5_4)\\
	\left(2-\lambda_m \right)u(Z^5_4)&=&u(Z^5_3)+u(Z^5_5)\\
	\left(2-\lambda_m \right)u(Z^5_5)&=&u(Z^5_4)+u(Z^5_6)\\
	\left(2-\lambda_m \right)u(Z^5_6)&=&u(Z^5_5)+u(Z^5_7)\\
	\left(2-\lambda_m \right)u(Z^5_7)&=&u(Z^5_6)+u(X_5)\\
	\end{array}\right. $$
	
	$$\left \lbrace \begin{array}{ccc}
	\left(2-\lambda_m \right)u(Z^6_1)&=&u(X_5)+u(Z^6_2)\\
	\left(2-\lambda_m \right)u(Z^6_2)&=&u(Z^6_1)+u(Z^6_3)\\
	\left(2-\lambda_m \right)u(Z^6_3)&=&u(Z^6_2)+u(Z^6_4)\\
	\left(2-\lambda_m \right)u(Z^6_4)&=&u(Z^6_3)+u(Z^6_5)\\
	\left(2-\lambda_m \right)u(Z^6_5)&=&u(Z^6_4)+u(Z^6_6)\\
	\left(2-\lambda_m \right)u(Z^6_6)&=&u(Z^6_5)+u(Z^6_7)\\
	\left(2-\lambda_m \right)u(Z^6_7)&=&u(Z^6_6)+u(X_6)\\
	\end{array}\right. \quad , \quad
	\left \lbrace \begin{array}{ccc}
	\left(2-\lambda_m \right)u(Z^7_1)&=&u(X_6)+u(Z^7_2)\\
	\left(2-\lambda_m \right)u(Z^7_2)&=&u(Z^7_1)+u(Z^7_3)\\
	\left(2-\lambda_m \right)u(Z^7_3)&=&u(Z^7_2)+u(Z^7_4)\\
	\left(2-\lambda_m \right)u(Z^7_4)&=&u(Z^7_3)+u(Z^7_5)\\
	\left(2-\lambda_m \right)u(Z^7_5)&=&u(Z^7_4)+u(Z^7_6)\\
	\left(2-\lambda_m \right)u(Z^7_6)&=&u(Z^7_5)+u(Z^7_7)\\
	\left(2-\lambda_m \right)u(Z^7_7)&=&u(Z^7_6)+u(X_7)\\
	\end{array}\right. $$
	
	$$\left \lbrace \begin{array}{ccc}
	\left(2-\lambda_m \right)u(Z^8_1)&=&u(X_7)+u(Z^8_2)\\
	\left(2-\lambda_m \right)u(Z^8_2)&=&u(Z^8_1)+u(Z^8_3)\\
	\left(2-\lambda_m \right)u(Z^8_3)&=&u(Z^8_2)+u(Z^8_4)\\
	\left(2-\lambda_m \right)u(Z^8_4)&=&u(Z^8_3)+u(Z^8_5)\\
	\left(2-\lambda_m \right)u(Z^8_5)&=&u(Z^8_4)+u(Z^8_6)\\
	\left(2-\lambda_m \right)u(Z^8_6)&=&u(Z^8_5)+u(Z^8_7)\\
	\left(2-\lambda_m \right)u(Z^8_7)&=&u(Z^8_6)+u(P_1)\\
	\end{array}\right. $$
	
	\noindent The system can be written as: $\mathbf{A}_2\mathbf{x}=\mathbf{0}$ and we look for the kernel of $\mathbf{A}_2$ for forbidden eigenvalues.\\
	
	\noindent $\rightsquigarrow$ \underline{For $\lambda_2=2$, the eigenspace is one dimensional, generated by the vector:}
	
	\begin{align*}
	&( 0 , 0 , 0 , 0 , 0 , 0 , 0 , -1 , 0 , 1 , 0 , -1 , 0 , 1 , -1 , 0 , 1 , 0 , -1 , 0 ,  1 , -1 , 0 , 1 , 0 , -1 , 0 , 1 , -1 , 0 , 1 , 0 , -1 , 0 , 1 , -1 , 0 , 1 , 0 ,\\
	&-1 , 0 , 1 , -1 , 0 , 1 , 0 , -1 , 0 , 1 , -1 , 0 , 1 , 0 , -1 , 0 , 1 , -1 , 0 ,
	1 , 0 , -1 , 0 , 1)
	\end{align*}
	
	\noindent $\rightsquigarrow$ \underline{For~$\lambda_2=2 +\varepsilon \sqrt{2}, \, \varepsilon\in \{-1,1\}$, the eigenspace is one dimensional, and is generated by:}
	
	\begin{align*}
	&(0 , 0 , 0 , 0 , 0 , 0 , 0 , -1 , \varepsilon\sqrt{2} , -1 , 0 , 1 , -\varepsilon\sqrt{2} , 1 , -1 ,
	\varepsilon\sqrt{2} , -1 , 0 , 1 , -\varepsilon\sqrt{2} , 1 , -1 , \varepsilon\sqrt{2} , -1 , 0 , 1 , -\varepsilon\sqrt{2} ,
	1 , -1 ,  \varepsilon\sqrt{2} , -1 ,  0 , 1 ,\\
	&  -\varepsilon\sqrt{2} , 1 , -1 , \varepsilon\sqrt{2} , -1 , 0 , 1 ,
	-\varepsilon\sqrt{2} , 1 , -1 , \varepsilon\sqrt{2} , -1 , 0 , 1 , -\varepsilon\sqrt{2} , 1 , -1 , \varepsilon\sqrt{2} , -1 ,
	0 , 1 , -\varepsilon\sqrt{2} , 1 , -1 , \varepsilon\sqrt{2} , -1 , 0 , 1 , -\varepsilon\sqrt{2} , 1)
	\end{align*}
	
	
	
	\noindent $\rightsquigarrow$ \underline{For~$\lambda_2=2 + \varepsilon\sqrt{2 + \sqrt{2}}, \, \varepsilon\in\{-1,1\}$, the eigenspace has dimension eight, and is generated by:}
	
	\begin{align*}
	&(0 , 0 , 0 , 0 , 0 , 0 , 0 , -1 , \varepsilon\sqrt{2+\sqrt{2}} , -1-\sqrt{2} , \varepsilon\sqrt{2
		\left(2+\sqrt{2}\right)} , -1-\sqrt{2} , \sqrt{2+\sqrt{2}} , -1 , 1 ,  -\varepsilon\sqrt{2+\sqrt{2}} , 1+\sqrt{2} ,\\
	& -\varepsilon\sqrt{2 \left(2+\sqrt{2}\right)} , 1+\sqrt{2}, -\sqrt{2+\sqrt{2}} , 1 , -1 , \varepsilon\sqrt{2+\sqrt{2}} , -1-\sqrt{2} , \varepsilon\sqrt{2
		\left(2+\sqrt{2}\right)} , -1-\sqrt{2} , \varepsilon\sqrt{2+\sqrt{2}} ,\\
	&   -1 , 1 ,
	-\varepsilon\sqrt{2+\sqrt{2}} , 1+\sqrt{2} , -\varepsilon\sqrt{2 \left(2+\sqrt{2}\right)} , 1+\sqrt{2}
	, -\varepsilon\sqrt{2+\sqrt{2}} , 1 , -1 , \varepsilon\sqrt{2+\sqrt{2}} , -1-\sqrt{2} , \varepsilon\sqrt{2  \left(2+\sqrt{2}\right)} ,\\
	& -1-\sqrt{2} , \varepsilon\sqrt{2+\sqrt{2}} , -1 , 1 ,
	-\varepsilon\sqrt{2+\sqrt{2}} , 1+\sqrt{2} , -\varepsilon\sqrt{2 \left(2+\sqrt{2}\right)} , 1+\sqrt{2}
	, -\varepsilon\sqrt{2+\sqrt{2}} , 1 , -1 , \varepsilon\sqrt{2+\sqrt{2}} ,\\
	& -1-\sqrt{2} , \varepsilon\sqrt{2
		\left(2+\sqrt{2}\right)} , -1-\sqrt{2} , \varepsilon\sqrt{2+\sqrt{2}} , -1 , 1 ,
	-\varepsilon\sqrt{2+\sqrt{2}} ,    1+\sqrt{2} , -\varepsilon\sqrt{2 \left(2+\sqrt{2}\right)} , 1+\sqrt{2}, \\
	& -\varepsilon\sqrt{2+\sqrt{2}} , 1) \\
	\end{align*}
	
	
	
	\noindent $\rightsquigarrow$ \underline{For~$\lambda_2=2 + \varepsilon\sqrt{2 - \sqrt{2}}, \, \varepsilon\in\{-1,1\}$, the eigenspace has dimension eight, and is generated by:}
	
	\begin{align*}
	&(0 , 0 , 0 , 0 , 0 , 0 , 0 , -1 , \varepsilon\sqrt{2-\sqrt{2}} , \sqrt{2}-1 , -\varepsilon\sqrt{2
		\left(2-\sqrt{2}\right)} , \sqrt{2}-1 , \varepsilon\sqrt{2-\sqrt{2}} , -1 , 1 ,  -\varepsilon\sqrt{2-\sqrt{2}} , 1-\sqrt{2} ,\\
	& \varepsilon\sqrt{2 \left(2-\sqrt{2}\right)} , 1-\sqrt{2} ,
	-\varepsilon\sqrt{2-\sqrt{2}} , 1 , -1 , \varepsilon\sqrt{2-\sqrt{2}} , \sqrt{2}-1 , -\varepsilon\sqrt{2
		\left(2-\sqrt{2}\right)} , \sqrt{2}-1 , \varepsilon\sqrt{2-\sqrt{2}} , -1 , 1 ,\\
	&   -\varepsilon\sqrt{2-\sqrt{2}} , 1-\sqrt{2} , \varepsilon\sqrt{2 \left(2-\sqrt{2}\right)} , 1-\sqrt{2} ,
	-\varepsilon\sqrt{2-\sqrt{2}} , 1 , -1 , \varepsilon\sqrt{2-\sqrt{2}} , \sqrt{2}-1 , -\varepsilon\sqrt{2
		\left(2-\sqrt{2}\right)} , \sqrt{2}-1 ,\\
	& \varepsilon\sqrt{2-\sqrt{2}} , -1 , 1 ,
	-\varepsilon\sqrt{2-\sqrt{2}} , 1-\sqrt{2} , \varepsilon\sqrt{2 \left(2-\sqrt{2}\right)} , 1-\sqrt{2} ,
	-\varepsilon\sqrt{2-\sqrt{2}} , 1 , -1 , \varepsilon\sqrt{2-\sqrt{2}} , \sqrt{2}-1 ,\\
	& -\varepsilon\sqrt{2 \left(2-\sqrt{2}\right)} , \sqrt{2}-1 , \varepsilon\sqrt{2-\sqrt{2}} , -1 , 1 ,
	-\varepsilon\sqrt{2-\sqrt{2}} , 1-\sqrt{2} , \varepsilon\sqrt{2 \left(2-\sqrt{2}\right)} , 1-\sqrt{2} ,
	-\varepsilon\sqrt{2-\sqrt{2}} , 1) \\
	\end{align*}
	
	
	
	\noindent From every forbidden eigenvalue~$\lambda_1$, the spectral decimation leads to eight eigenvalues. Each of these eigenvalues has multiplicity~$1$.\\
	
	\noindent One may easily check that:

	$$\# V_2\setminus V_0 =63= 7\times 1 + 8\times 7$$
	
	\noindent Thus, the spectrum is complete.

	\vskip 1cm

	\subsubsection{General case}

	\noindent Let us now go back to the general case. Given a strictly positive integer~$m$, let us introduce the respective multiplicities~$M_m(\lambda_m)$ of the eigenvalue $\lambda_m$.\\
	
	\noindent One can easily check by induction that:
	
	$$\# V_m\setminus V_0 =8^m -1$$

	\begin{theorem}
		To every forbidden eigenvalue~$\lambda$ is associated an eigenspace, the dimension of which is one: $M_m(\lambda)=1$.
	\end{theorem}

	\vskip 1cm
	\begin{proof}
		The theorem is true for $m=1$ and $m=2$.
		Recursively, we suppose the result true for $m-1$.\\
		At $m$ There are :
		$$\sum_{i=1}^{m-1} 8^i \times 7=8^{m} - 8=\# V_m\setminus V_0 - 7$$
		eigenvalues generated  using the continuous formula of spectral decimation. The remaining eigenspace has dimension $7$.\\
		We can verify that every forbidden eigenvalue $\lambda$ is an eigenvalue for $m$. If we consider the eigenfunction null everywhere except the points of $V_m\setminus V_{m-1}$ where it takes the values of $V^1_{\lambda}$ in the interior of every $m$-cell.\\
		We conclude that the eigenspace of every forbidden eigenvalue is one dimensional.
	\end{proof}
	
	\begin{figure}[h!]
		\center{\psfig{height=7cm,width=7cm,angle=0,file=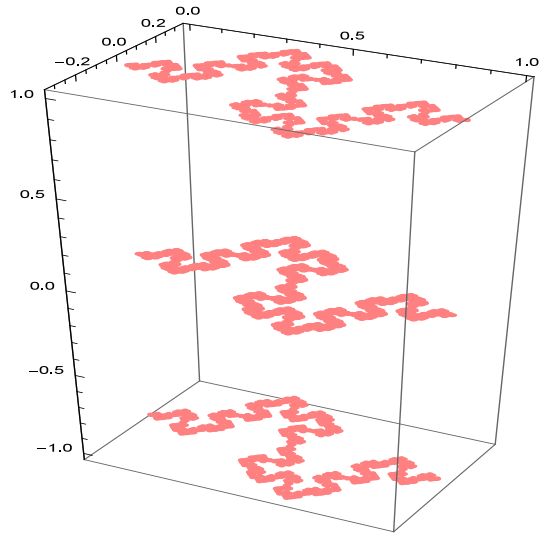}}
		\caption{The graph of the first eigenfunction}
	\end{figure}
	
	\begin{figure}[h!]
		\center{\psfig{height=7cm,width=7cm,angle=0,file=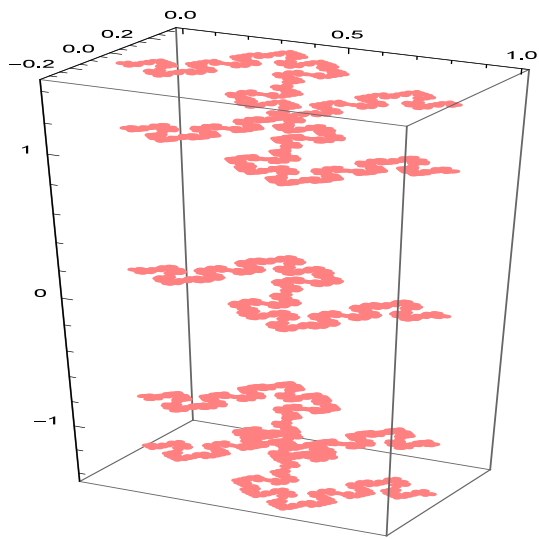}}
		\caption{The graph of the second eigenfunction}
	\end{figure}
	
	\begin{figure}[h!]
		\center{\psfig{height=7cm,width=7cm,angle=0,file=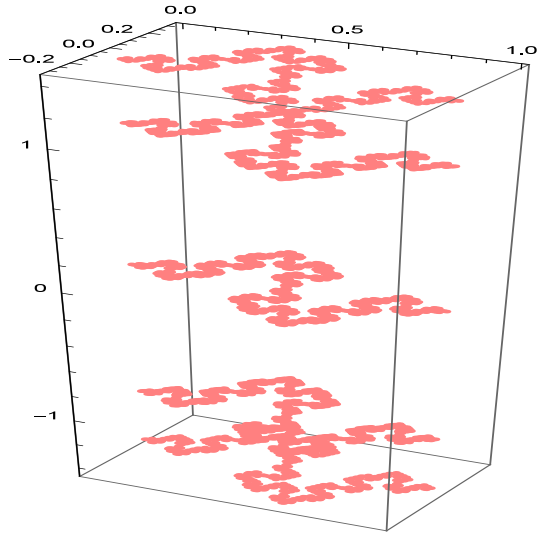}}
		\caption{The graph of the third eigenfunction}
	\end{figure}
	
	\begin{figure}[h!]
		\center{\psfig{height=7cm,width=7cm,angle=0,file=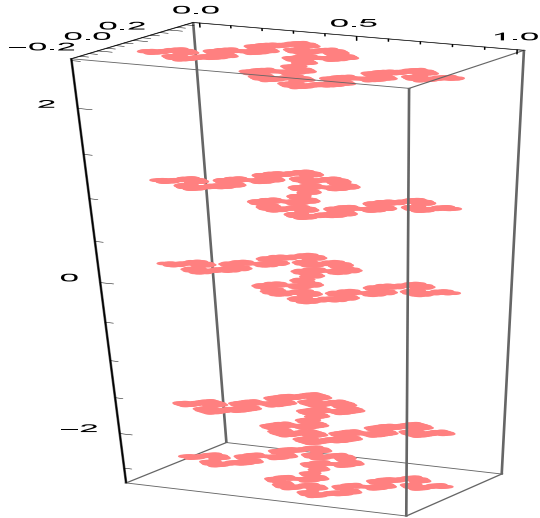}}
		\caption{The graph of the fourth eigenfunction}
	\end{figure}
	
	\begin{figure}[h!]
		\center{\psfig{height=7cm,width=7cm,angle=0,file=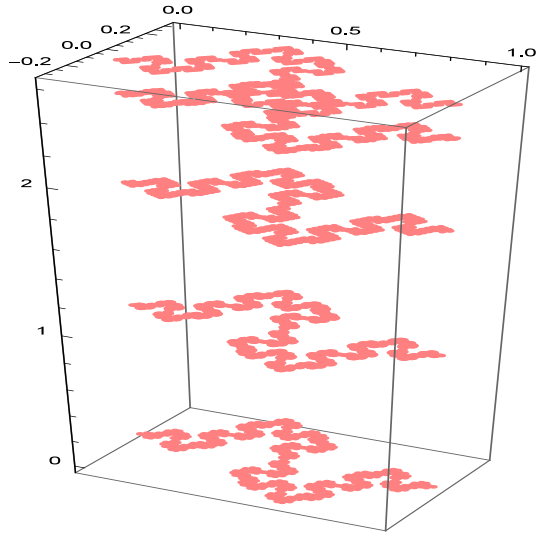}}
		\caption{The graph of the fifth eigenfunction}
	\end{figure}
	
	\begin{figure}[h!]
		\center{\psfig{height=7cm,width=7cm,angle=0,file=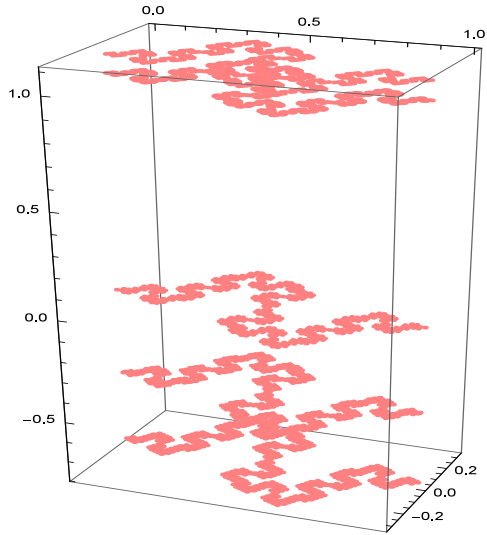}}
		\caption{The graph of the sixth eigenfunction}
	\end{figure}
	
	\begin{figure}[h!]
		\center{\psfig{height=7cm,width=7cm,angle=0,file=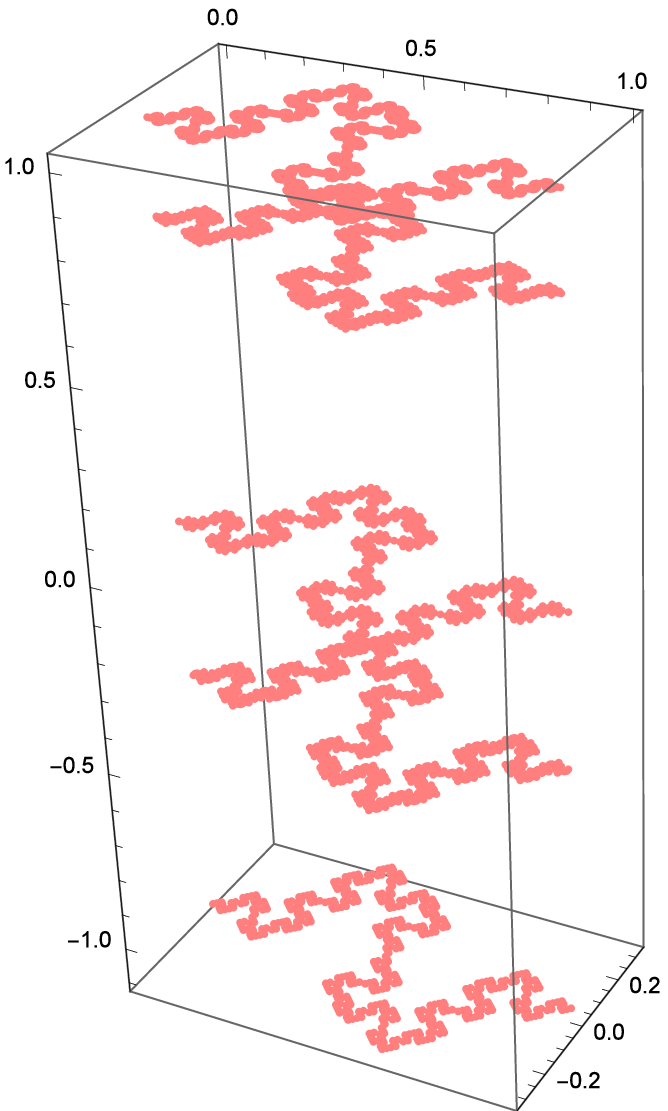}}
		\caption{The graph of the seventh eigenfunction}
	\end{figure}
	
	\newpage
	\section{Metric - Towards spectral asymptotics}
	
	\begin{definition}\textbf{Effective resistance metric, on~${\mathfrak {MC}}$}\\
		\noindent Given two points~$(X,Y)$ of~\mbox{${\mathfrak {MC}}^2$}, let us introduce the \textbf{effective resistance metric between~$X$ and~$Y$}:
		
		$$ R_{\mathfrak {MC}}(X,Y)= \left \lbrace \min_{\left \lbrace u \, |\, u(X)=0 , u(Y)=1\right \rbrace } \,{\cal E}(u )  \right \rbrace^{-1}$$
		
		\noindent In an equivalent way,~\mbox{$ R_{\mathfrak {MC}}(X,Y)$} can be defined as the minimum value of the real numbers~$R$ such that, for any function~$u$ of~\mbox{$\text{dom} \,\Delta$}:

		$$ \left |u (X)-u (Y)\right |^2 \leq R\,  {\cal E}(u)  $$
		
	\end{definition}

	\vskip 1cm
	
	\begin{definition}\textbf{Metric, on the Minkowski Curve~${\mathfrak {MC}}$}\\
		\noindent Let us define, on the Minkowski Curve~${\mathfrak {MC}}$, the distance~\mbox{$d_{ \mathfrak {MC} }$} such that, for any pair of points~$(X,Y)$
		of~\mbox{$ {\mathfrak {MC}}^2 $}:
		
		$$ d_{{\mathfrak {MC}}}(X,Y)= \left \lbrace \min_{\left \lbrace u \, |\, u(X)=0 , u(Y)=1\right \rbrace } \,{\cal E}(u,u)  \right \rbrace^{-1}$$
	\end{definition}
	
	\vskip 1cm

	\begin{remark}
		
		One may note that the minimum
		
		$$\min_{\left \lbrace u \, |\, u(X)=0 , u(Y)=1\right \rbrace } \,{\cal E}(u )   $$
		
		\noindent is reached for~$u$ being harmonic on the complement set, on~\mbox{$\mathfrak {MC}$}, of the set
		$$ \left \lbrace X \right \rbrace \cup  \left \lbrace Y \right \rbrace  $$
		
		\noindent (One might bear in mind that, due to its definition, a harmonic function~$u$ on~\mbox{$\mathfrak {MC}$} minimizes the sequence of energies~\mbox{$\left ({\cal E}_{ {\mathfrak {MC}}_m }  (u,u )\right)_{m\in\N}$}.\\
		
	\end{remark}
	
	\vskip 1cm

	\begin{definition}\textbf{Dimension of the Minkowski Curve~${\mathfrak {MC}}$, in the resistance metrics}\\
		
		\noindent The \textbf{dimension of the Minkowski Curve~${\mathfrak {MC}}$}, in the resistance metrics, is the strictly positive number~\mbox{$d_{\mathfrak {MC} }$} such that, given a strictly positive real number~$r$, and a point~$X\,\in\,{\mathfrak {MC}} $, for the~\mbox{$X-$centered} ball of radius~$r$, denoted by~\mbox{${\cal B}_r(X)$}:
		
		$$ \mu \left ({\cal B}_r(X) \right) =r^{d_{\mathfrak {MC}}}$$
		
	\end{definition}
	
	\vskip 1cm
	
	\begin{pte}
		
		\noindent Given a natural integer~$m$, and two points~$(X,Y)$ of~\mbox{$ {\mathfrak {MC}}^2 $} such that~\mbox{$X  \underset{m }{\sim}  Y $}:
		
		$$ \displaystyle \min_{\left \lbrace u \, |\, u(X)=0 , u(Y)=1\right \rbrace } \,{\cal E}(u ) \lesssim r^m\, 4^{-2\, m \,\delta} =  \displaystyle \frac{1}{8^m}  $$

		
		


		\noindent Let us denote by~$\mu$ the standard measure on~$\mathfrak {MC}$ which assigns measure~\mbox{$\displaystyle \frac{1}{4^{2m}}$ } to each quadrilateral~$m-$ cell. Let us now look for a real number~$d_{\mathfrak {MC}}$ such that:

		$$\left (\displaystyle \frac{1}{8 }\right)^{m \,d_{\mathfrak {MC}}}= \displaystyle \frac{1}{4^{2m}}$$

		\noindent One obtains:

		$$d_{\mathfrak {MC} }=\displaystyle \frac{2}{3}$$

		\noindent   Given a strictly positive real number~$r$, and a point~$X\,\in\,{\mathfrak {MC}} $, one has then the following estimate, for the~\mbox{$X-$centered} ball of radius~$r$, denoted by~\mbox{${\cal B}_r(X)$}:
		
		$$ \mu \left ({\cal B}_r(X) \right) =r^{d_{\mathfrak {MC}}}$$
		
	\end{pte}

	\vskip 1cm

	\begin{definition}\textbf{Eigenvalue counting function}\\
		
		\noindent We introduce the eigenvalue counting function~${\cal N}^{\mathfrak {MC}}$ such that, for any real number~$x$:
		
		$${\cal N}^{\mathfrak {MC}}(x) = \#\, \left \lbrace \lambda\,\text{eigenvalue of $-\Delta$} \, : \quad \lambda \leq  x \right \rbrace $$

	\end{definition}

	\vskip 1cm

	\begin{pte}
		\noindent According to J.~Kigami~\cite{KigamiWeylFormula}, one has the modified Weyl formula:

		$${\cal N}^{\mathfrak {MI}}(x)= \left( G\left (\frac{\ln x}{2}\right )+o(1)\right) {x^{D_{S}\left(\mathfrak{MI}\right)}}$$

	\end{pte}
	
	\vskip 1cm
	
	\section{From the Minkowski Curve, to the Minskowski Island}
	
	By connecting four Minkowski Curves, as it can be seen on the following figure:

	\begin{figure}[h!]
		\center{\psfig{height=6cm,width=6cm,angle=0,file=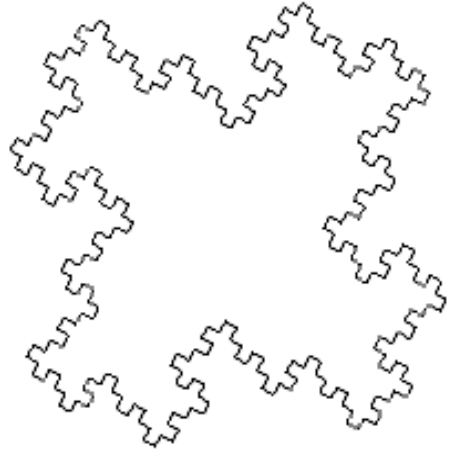}}
		\caption{The Minskowski Island.}
	\end{figure}

	\noindent one obtains what is called "the Minskowski Island".\\
	
	\noindent We expose, in the sequel, our results (we refer to~Zhou~\cite{Zhou2007} for further details).\\
	In the case of the Minkowski Curve, one may check that:
	
	\[ \#(F_i(V_0)\cap V_0)\leq 1 \quad \forall \ i \in \{1,2,3,4,5,6,7,8\} \]

	$$
	D=
	\left(
	\begin{array}{cccc}
	-2 & 1 & 0 & 1\\
	1 & -2 & 1 & 0\\
	0 & 1 & -2 & 1 \\
	1 & 0 & 1 & -2 \\
	\end{array}
	\right)
	\quad , \quad
	T=
	\left(
	\begin{array}{cccc}
	-2 & 0 & 0 & 0\\
	0 & -2 & 0 & 0\\
	0 & 0 & -2 & 0 \\
	0 & 0 & 0 & -2 \\
	\end{array}
	\right)
	$$
	Let us set:
	\[
	B=
	\left(
	\begin{array}{ccccccc}
	-2 & 1 & 0 & 0 & 0 & 0 & 0 \\
	1 & -2 & 1 & 0 & 0 & 0 & 0 \\
	0 & 1 & -2 & 1 & 0 & 0 & 0 \\
	0 & 0 & 1 & -2 & 1 & 0 & 0 \\
	0 & 0 & 0 & 1 & -2 & 1 & 0 \\
	0 & 0 & 0 & 0 & 1 & -2 & 1 \\
	0 & 0 & 0 & 0 & 0 & 1 & -2 \\
	\end{array}
	\right)
	\]
	and:
	
	\[
	C=
	\left(
	\begin{array}{ccccccc}
	2 & 0 & 0 & 0 & 0 & 0 & 0 \\
	0 & 2 & 0 & 0 & 0 & 0 & 0 \\
	0 & 0 & 2 & 0 & 0 & 0 & 0 \\
	0 & 0 & 0 & 2 & 0 & 0 & 0 \\
	0 & 0 & 0 & 0 & 2 & 0 & 0 \\
	0 & 0 & 0 & 0 & 0 & 2 & 0 \\
	0 & 0 & 0 & 0 & 0 & 0 & 2 \\
	\end{array}
	\right)
	\]
	
	Then:
	
	$$
	X=
	\left(
	\begin{array}{cccc}
	B & 0 & 0 & 0 \\
	0 & B & 0 & 0 \\
	0 & 0 & B & 0 \\
	0 & 0 & 0 & B \\
	\end{array}
	\right)
	$$
	
	and:
	
	$$
	M=
	\left(
	\begin{array}{cccc}
	C & 0 & 0 & 0 \\
	0 & C & 0 & 0 \\
	0 & 0 & C & 0 \\
	0 & 0 & 0 & C \\
	\end{array}
	\right)
	\quad , \quad
	J=
	\left(
	\begin{array}{cccc}
	1 & 0 & 0 & 0 \\
	0 & 0 & 0 & 0 \\
	0 & 0 & 0 & 0 \\
	0 & 0 & 0 & 0 \\
	0 & 0 & 0 & 0 \\
	0 & 0 & 0 & 0 \\
	0 & 1 & 0 & 0 \\
	0 & 1 & 0 & 0 \\
	0 & 0 & 0 & 0 \\
	0 & 0 & 0 & 0 \\
	0 & 0 & 0 & 0 \\
	0 & 0 & 0 & 0 \\
	0 & 0 & 0 & 0 \\
	0 & 0 & 1 & 0 \\
	0 & 0 & 1 & 0 \\
	0 & 0 & 0 & 0 \\
	0 & 0 & 0 & 0 \\
	0 & 0 & 0 & 0 \\
	0 & 0 & 0 & 0 \\
	0 & 0 & 0 & 0 \\
	0 & 0 & 0 & 1 \\
	0 & 0 & 0 & 1 \\
	0 & 0 & 0 & 0 \\
	0 & 0 & 0 & 0 \\
	0 & 0 & 0 & 0 \\
	0 & 0 & 0 & 0 \\
	0 & 0 & 0 & 0 \\
	1 & 0 & 0 & 0 \\
	\end{array}
	\right)
	$$
	
	\noindent One may check that the Minkowsky Island has a strong harmonic structure, and that:
	\footnotesize
	\[T + J^T \,(X +\lambda  M)^{-1}J=-\frac{1}{8
		(\lambda -1)\, (2 \,(\lambda -2)\, \lambda +1) \left(8 \,(\lambda -2) \lambda \, (\lambda -1)^2+1\right)}\,D+
	\frac{\lambda \, (4 \,\lambda  (\lambda  (2 \lambda -7)+7)-7)}{8\, (\lambda -2) \lambda \, (\lambda -1)^2+1}\,T\]
	\normalsize
	\noindent We obtain, as previously:
	
	\[\frac{1}{K_D(0)}=8\]
	
	\noindent The spectral decimation function is given, for any real number~$\lambda$, by:
	
	\[
	R(\lambda)  =-(\lambda-4) (\lambda-2)^2 \lambda ((\lambda-4) \lambda + 2)^2
	\]
	
	\noindent which leads to the same spectrum as the one of the Curve.
	
	
	\vskip 1cm

	\bibliographystyle{alpha}
	\bibliography{BibliographieClaire}

\end{document}